\documentclass[a4,10pt]{article}
\usepackage{amsmath,amscd,amssymb}

\newcommand{\nbiga}{\mathcal{A}}

\newcommand{\nbigd}{\mathcal{D}}
\newcommand{\nbige}{\mathcal{E}}
\newcommand{\nbigf}{\mathcal{F}}
\newcommand{\nbigg}{\mathcal{G}}

\newcommand{\nbigk}{\mathcal{K}}

\newcommand{\nbigo}{\mathcal{O}}
\newcommand{\nbigp}{\mathcal{P}}

\newcommand{\nbigu}{\mathcal{U}}

\newcommand{\seisuu}{\mathbb{Z}}

\newcommand{\cnum}{{\boldsymbol C}}

\newcommand{\EE}{\mathbb{E}}

\newcommand{\gbigk}{\mathfrak K}

\newcommand{\gminia}{\mathfrak a}
\newcommand{\gminib}{\mathfrak b}

\newcommand{\gminik}{\mathfrak k}

\newcommand{\vecv}{{\boldsymbol v}}

\newcommand{\vecf}{{\boldsymbol f}}

\newcommand{\vecg}{{\boldsymbol g}}

\newcommand{\lrarr}{\longrightarrow}

\newcommand{\pf}{{\bf Proof}\hspace{.1in}}
\newcommand{\qed}{\mbox{\rule{1.2mm}{3mm}}}

\def\End{\mathop{\rm End}\nolimits}

\def\GL{\mathop{\rm GL}\nolimits}

\def\rank{\mathop{\rm rank}\nolimits}
\def\Spec{\mathop{\rm Spec}\nolimits}

\def\modulo{\mathop{\rm mod}\nolimits}

\def\Res{\mathop{\rm Res}\nolimits}

\def\ord{\mathop{\rm ord}\nolimits}

\def\Tr{\mathop{\rm Tr}\nolimits}

\def\id{\mathop{\rm id}\nolimits}

\def\Irr{\mathop{\rm Irr}\nolimits}

\newcommand{\del}{\partial}

\newcommand{\etabar}{\overline{\eta}}
\newcommand{\nablabar}{\overline{\nabla}}

\newcommand{\laplacianlambda}

\newcommand{\Sp}{{\mathcal Sp}}

\def\reg{\mathop{\rm reg}\nolimits}

\def\disc{\mathop{\rm disc}\nolimits}
\def\Fr{\mathop{\rm Fr}\nolimits}
\def\Sol{\mathop{\rm Sol}\nolimits}

\newcommand{\nablahat}{\widehat{\nabla}}

\newcommand{\kappatilde}{\widetilde{\kappa}}

\newcommand{\Dhat}{\widehat{D}}

\newcommand{\kuebar}{\overline{k}}
\newcommand{\Xtilde}{\widetilde{X}}
\newcommand{\Utilde}{\widetilde{U}}
\newcommand{\huebar}{\overline{h}}
\newcommand{\gminiatilde}{\widetilde{\gminia}}
\newcommand{\Xbar}{\overline{X}}

\newcommand{\nbigutilde}{\widetilde{\nbigu}}
\newcommand{\nbigebar}{\overline{\nbige}}
\newcommand{\nbigehat}{\widehat{\nbige}}
\newcommand{\Ptilde}{\widetilde{P}}

\newcommand{\Fpbar}{\overline{\mathbb F}_p}
\newcommand{\boldnbigp}{\boldsymbol\nbigp}
\newcommand{\Phat}{\widehat{P}}
\newcommand{\Atilde}{\widetilde{A}}
\newcommand{\Ytilde}{\widetilde{Y}}
\newtheorem{thm}{Theorem}[section]
\newtheorem{cor}[thm]{Corollary}

\newtheorem{lem}[thm]{Lemma}
\newtheorem{prop}[thm]{Proposition}
\newtheorem{df}[thm]{Definition}

\newtheorem{assumption}[thm]{Assumption}

\begin{document}

\title{Good formal structure for
 meromorphic flat connections\\
 on smooth projective surfaces}
\author{Takuro Mochizuki}
\date{}
\maketitle

\begin{abstract}

We prove the algebraic version
of a conjecture of C. Sabbah
on the existence of the good formal structure
for meromorphic flat connections on surfaces
after some blow up.

\vspace{.1in}

\end{abstract}

\section{Introduction}

\subsection{Main result}

Let $X$ be a smooth complex projective surface,
and let $D$ be a normal crossing divisor of $X$.
Let $(\nbige,\nabla)$ be
a flat meromorphic connection 
on $(X,D)$,
i.e.,
$\nbige$ denotes a locally free $\nbigo_X(\ast D)$-module,
and $\nabla:\nbige\lrarr\nbige\otimes\Omega_{X/\cnum}^1$
denotes a flat connection.
We discuss a conjecture of C. Sabbah
under the algebraicity assumption.

\begin{thm}
\label{thm;07.4.15.10}
There exists a regular birational morphism
$\pi:\Xtilde\lrarr X$ such that 
$\pi^{-1}(\nbige,\nabla)$ has the good formal structure.
\end{thm}

See Subsection \ref{subsection;07.4.30.1}
for good formal structure.
For explanation of the meaning of the theorem,
let us recall the very classical result in the curve case.
(See the introduction of \cite{majima} for more detail,
for example.)
Let $C$ be a smooth projective curve,
and let $Z\subset C$ be a finite subset.
Let $(\nbige,\nabla)$ be a meromorphic connection on $(C,Z)$,
i.e., $\nbige$ is a locally free $\nbigo(\ast Z)$-module
with a connection $\nabla$.
Let $P$ be any point of $Z$,
and let $(U,t)$ be a holomorphic coordinate neighbourhood
around $P$ such that $t(P)=0$.
The local structure of $\nbige$ around $P$
can be understood by
the formal structure at $P$
and the Stokes structure around $P$.
Namely,
take a ramified covering
$\varphi_P:(\Utilde,t_d)\lrarr (U,t)$ given by $t=t_d^d$,
where $d$ is a large integer 
divided by $(\rank(\nbige)!)^3$, for example.
Let $\Ptilde\in\Utilde$ be the inverse image of $P$.
The formal completion of $\varphi_P^{\ast}(\nbige,\nabla)$
at $\Ptilde$ is decomposed into the direct sum
$\bigoplus_{\gminia\in \cnum(\!(t_d)\!)/\cnum[\![t_d]\!]}
 (\nbige_{\gminia},\nabla_{\gminia})$,
where $\nabla_{\gminia}-d\gminia$ are regular.
(In the curve case,
we do not have to assume that the base space is
algebraic.
In fact, the decomposition can be obtained
for any connections on formal curves.)
Then, the formal decomposition
can be lifted to the decomposition
on any small sectors by the asymptotic analysis,
which leads us the Stokes structure.

It is a challenging and foundational problem to obtain the generalization
in the higher dimensional case.
The systematic study 
was initiated by H. Majima,
who developed the asymptotic analysis in the higher dimensional case.
(See \cite{majima}, for example.)
Briefly speaking,
his result gives the lifting of a formal decomposition
to the decomposition on small sectors.
Inspired by Majima's work,
Sabbah (\cite{sabbah4}) developed the asymptotic analysis
in the other framework.
He observed the significance of the understanding
on the formal structure of the irregular connection.
He proposed the conjecture which says
that Theorem \ref{thm;07.4.15.10}
may hold without the algebraicity assumption,
and he established it in the case $\rank(\nbige)\leq 5$.
He also reduced the problem to the study
of the turning points contained in the smooth part
of the divisor $D$,
without any assumption on the rank.

Sabbah gave some interesting applications of the conjecture,
one of which is a conjecture of B. Malgrange
on the absence of the confluence phenomena for 
flat meromorphic connections.
Recently, Y. Andr\'e (\cite{andre}) proved Malgrange's conjecture
motivated by Sabbah's conjecture.

In this paper,
we will give a proof of the algebraic version of Sabbah's conjecture.
In \cite{mochi6},
the author intends to establish the correspondence
of semisimple algebraic holonomic $D$-modules
and polarizable wild pure twistor $D$-modules
through wild harmonic bundles,
on smooth projective surfaces and
the higher dimensional varieties,
which is related with a conjecture of M. Kashiwara \cite{k5}.
Theorem \ref{thm;07.4.15.10} has the foundational importance
for the study.

\subsection{Main ideas}

Let $k$ be an algebraically closed field,
and let $(\nbige,\nabla)$ be a meromorphic flat connection on 
$k[\![s]\!](\!(t)\!)$.
If the characteristic number $p$ of $k$ is positive,
we always assume that $p$ is much larger than
$\rank\nbige$ and the Poincar\'e rank of $\nbige$
with respect to $t$.
Let $\nabla_t$ denote 
the induced relative connection
$\nbige\lrarr\nbige\otimes
 \Omega^1_{k[\![s]\!](\!(t)\!)/k[\![s]\!]}$.
The induced connection
$(\nbige,\nabla_t)\otimes k(\!(s)\!)(\!(t)\!)$
is denoted by $(\nbige_1,\nabla_1)$.
The specialization of $(\nbige,\nabla_t)$
at $s=0$ is denoted by $(\nbige_0,\nabla_0)$.
We have the set of the irregular values
$\Irr(\nbige_0,\nabla_0)
\subset k(\!(t_d)\!)/k[\![t_d]\!]$
and 
$\Irr(\nbige_1,\nabla_1)
\subset
 k(\!(s_d)\!)(\!(t_d)\!)/k(\!(s_d)\!)[\![t_d]\!]$,
where $t_d$ and $s_d$ denote $d$-th roots of $t$ and $s$,
respectively.
Briefly and imprecisely speaking,
one of the main issues is how to compare
$\Irr(\nbige_i,\nabla_i)$ $(i=0,1)$.
Ideally, we hope that
$\Irr(\nbige_1,\nabla_1)$
is contained in $k[\![s]\!](\!(t_d)\!)/k[\![s,t_d]\!]$,
and that the specialization at $s=0$ gives
$\Irr(\nbige_0,\nabla_0)$.
However, they are not true, in general.

In the case $p>0$,
we have the $p$-curvature $\psi$ (resp. $\psi_i$) of 
the connection $\nabla$
(resp. $\nabla_i$).
\[
 \psi\in \End(\nbige)\otimes F^{\ast}\Omega^1_{k[\![s]\!](\!(t)\!)/k},
\quad
 \psi_1\in\End(\nbige_1)\otimes
 F^{\ast}\Omega^1_{k(\!(s)\!)(\!(t)\!)/k(\!(s)\!)},
\quad
 \psi_0\in\End(\nbige_0)\otimes
 F^{\ast}\Omega^1_{k(\!(t)\!)/k}
\]
Here, $F$ denotes the absolute Frobenius map.
In the following,
we use the notation $\psi(t\del_t)$
to denote $\psi(F^{\ast}t\del_t)$,
for simplicity.
Let $\Sp(\psi(t\del_t))$ denote the set of the eigenvalues of
$\psi(t\del_t)$,
which is contained in $\nbiga_t$,
where 
$\nbiga$ denotes a finite extension of $k[\![s,t]\!]$ 
and $\nbiga_t$ denotes a localization of $\nbiga$
with respect to $t$.
Similarly,
let $\Sp(\psi_i(t\del_t))$ denote the set of the eigenvalues
of $\psi_i(t\del_t)$ for $i=0,1$, and then
$\Sp(\psi_0(t\del_t))\subset
 k(\!(t_d)\!)$ 
and $\Sp(\psi_1(t\del_t))\subset k(\!(s_d)\!)(\!(t_d)\!)$
for some appropriate $d\in\seisuu_{>0}$.
We may have the natural inclusion
$\kappa_1:\nbiga_t\lrarr k(\!(s_d)\!)(\!(t_d)\!)$
and the specialization
$\kappa_0:\nbiga_t\lrarr k(\!(t_d)\!)$ at $s=0$.
Clearly, $\psi_i(t\del_t)$ $(i=0,1)$ are naturally obtained from
$\psi(t\del_t)$ by $\kappa_i$,
and hence 
$\Sp(\psi_i(t\del_t))$ are obtained from
$\Sp(\psi(t\del_t))$ by $\kappa_i$.
Recall that the irregular value of $\nabla_i$ can be related
with the negative part of the eigenvalues
of $\psi_i(t\del_t)$ (Lemma \ref{lem;07.4.14.3}),
where the negative part of $f=\sum f_j\cdot t_d^j\in R(\!(t_d)\!)$
is defined to be $f_-:=\sum_{j<0}f_j\cdot t_d^j$.
Hence, we have the following diagram:
\[
 \begin{CD}
 \Sp\bigl(\psi_1(t\del_t)\bigr)
 @<{\kappa_1}<<
 \Sp\bigl(\psi(t\del_t)\bigr)
 @>{\kappa_0}>>
 \Sp\bigl(\psi_0(t\del_t)\bigr)\\
 @VVV @. @VVV\\
 \Irr(\nbige_1,\nabla_1)
@. @.
 \Irr(\nbige_0,\nabla_0)
 \end{CD}
\]
But, we should remark that
$\kappa_0(\alpha)_-$ and $\kappa_1(\alpha)_-$
cannot be directly related, in general.

Let us consider the simplest case where
the ramification of $\nbiga$ over $k[\![s,t]\!]$
may occur only at the divisor $\{t=0\}$.
Then, $\Sp(\psi(t\del_t))$ is contained in
$k[\![s]\!](\!(t_d)\!)$,
and $\kappa_0(\alpha)_-$ is the specialization of
$\kappa_1(\alpha)_-$ at $s=0$
for any $\alpha\in\Sp(\psi(t\del_t))$.
Thus, we can compare the irregular values
of $\nabla_i$ $(i=0,1)$ in this simplest case.

Then, we have to consider what happens
if the ramification of $\nbiga$ may be non-trivial.
As the second simplest case,
we assume that
the ramification may occur only
at the normal crossing divisor 
$(t)\cup(s')$ of $\Spec^fk[\![s,t]\!]$,
where $s'=s+t\cdot h(t)$.
Then, 
$\Sp(\psi(t\del_t))$ 
are contained in $k[\![s_d']\!](\!(t_d)\!)$,
where $s_d'$ denotes a $d$-th root of $s'$.
We assume, moreover, that
$\Sp(\psi(t\del_t))$ 
are contained in $k[\![s]\!](\!(t_d)\!)+k[\![s_d',t_d]\!]$.
Then, the negative part of the eigenvalues behave well
with respect to the specialization,
i.e.,
$\kappa_1(\alpha)_-=\kappa_0(\alpha)_-$
for any $\alpha\in \Sp(\psi(t\del_t))$.
Hence, we can compare the irregular values
of $\nabla_i$ $(i=0,1)$
in this mildly ramified case
(Lemma \ref{lem;07.4.14.10}).

We would like to apply such consideration to our problem.
Essentially, the problem is the following,
although we will discuss it in a different way.
Let $(\nbige,\nabla)$ be a meromorphic flat connection
with a lattice $E$ on $(X,D)$.
For simplicity, we assume 
everything is defined over $\seisuu$.
Then, we have the mod $p$-reductions
$(\nbige_p,\nabla_p):=(\nbige,\nabla)\otimes \Fpbar$
over $(X_p,D_p):=(X,D)\otimes\Fpbar$
with the lattice $E_p=E\otimes\Fpbar$.
Let $\psi_p\in\End(E_p)\otimes F^{\ast}\Omega_{X_p}(ND_p)$
denote the $p$-curvature.
Then, we have the spectral manifold
$ \Sigma_p(\psi_p):=
 \bigl\{(x,\omega)\,\big|\,
 \mbox{$\omega$ eigenvalues of $\psi_{p|x}$}
 \bigr\}
\subset 
 F^{\ast}\bigl(\Omega^1_{X_p}\otimes\nbigo(ND_p)\bigr)$.

For simplicity, we assume that
$\psi_p$ has the distinct eigenvalues
at the generic point.
Then, we hope that 
the ramification of the projection $\pi_p$
of $\Sigma_p(\psi_p)$ to $X_p$
may happen at normal crossing divisor,
after some blow ups,
i.e.,
$ R(\pi_p):=
 \bigl\{x\in X_p\,\big|\,\mbox{\rm $\pi_p$ is not etale at $x$}\bigr\}$
is normal crossing.
If we fix $p$,
it is easy to obtain such birational map
because we are considering the surface case.
But, for our problem,
we would like to control the ramification 
for almost all $p$ at once.
So we need something more.

Here, we recall the important observation of
J. Bost, Y. Laszlo and C. Pauly \cite{laszlo-pauly}
which says that 
we have $\Sigma_p'$ contained in
$\Omega^1_{X_p}\otimes\nbigo(ND)$,
such that
$\Sigma_p(\psi_p)$ is the pull back of $\Sigma_p'$.
So, we have only to control
the ramification curves $R(\pi_p')$
of the projection $\pi_p'$ of $\Sigma_p'$ to $X_{p}$.
Then, it is not difficult to see that
the arithmetic genus of $R(\pi_p')$ are dominated,
independently of $p$.
So, the complexity of the singularities
of these ramification curves are bounded,
and thus we can control them uniformly.
(See Section \ref{section;07.4.28.2}.)

\subsection{Acknowledgement}

The author thanks C. Sabbah for the discussion
and his attractive conjecture.
It is a pleasure to express his gratitude to
R. Bezrukavnikov,
O. Biquard,
M. Kontsevich, J. Li,
C. Simpson, A. Usnich,
D. Wei for some discussions.
He is grateful to Y. Tsuchimoto and A. Ishii
for their constant encouragement.
Special thanks goes to K. Vilonen.
He thanks the colleague in Kyoto University
for their support.
This paper is prepared during his visit at IHES and ICTP.
The author is grateful for their excellent hospitality.
He also thanks the partial financial supports by
Sasakawa Foundation
and Ministry of Education, Culture, Sports, Science and Technology.

The contents of this paper is 
a detailed proof of a theorem given 
in a talk at the conference
`Algebraic Analysis and Around' 
in honour of Professor Masaki Kashiwara's 60's birthday.
The author would like to express 
his gratitude to the organizers.

It is an extremely great pleasure 
for the author to dedicate this paper
to Masaki Kashiwara with admiration 
for his great works and 
his leading role in the development of
current mathematics.

\section{Preliminary}

\subsection{Notation}
\label{subsection;07.4.26.1}

Let $R$ be a ring, and let $t$ be a formal variable.
We use the notation $R[\![t]\!]$ (resp. $R(\!(t)\!)$)
to denote the ring of formal power series
(resp. the ring of formal Laurent power series)
over $R$.
Let $R(\!(t)\!)_{<0}$ denote the subset
$\bigl\{\sum_{j<0} a_j \cdot t^j\in R(\!(t)\!)\bigr\}$.
For any $f=\sum a_j\cdot t^j\in R(\!(t)\!)$,
we put $\ord_t(f):=\min\{j\,|\,a_j\neq 0\}$.
If we are given two variables $s$ and $t$,
we use the notation $R[\![s]\!](\!(t)\!)$
to denote the ring of formal Laurent power series
over $R[\![s]\!]$.
The notation $R(\!(t)\!)[\![s]\!]$ is used to denote
the ring of formal power series over $R(\!(t)\!)$.
We have
$R[\![s]\!](\!(t)\!)\subsetneq R(\!(t)\!)[\![s]\!]$.

For a given integer $d>0$ and a formal variable $t$,
we use the notation $t_d$ as a $d$-th root of $t$,
i.e., $t_d^d=t$.
For any $f=\sum f_j\cdot t_d^j\in R(\!(t_d)\!)$,
we put $f_-:=\sum_{j<0} f_j\cdot t_d^j$,
which is called the negative part of $f$.
If $d'$ is a factor of $d$,
we regard $R(\!(t_{d'})\!)$ as the subring of
$R(\!(t_{d})\!)$.
For any $f\in R(\!(t_d)\!)$,
we put $\ord_t(f):=d^{-1}\cdot\ord_{t_d}(f)$.
The definition is consistent for
the inclusions
$R(\!(t)\!)\subset R(\!(t_{d'})\!)\subset R(\!(t_{d})\!)$.
Let us consider the case where
$R$ is a ring over $\seisuu/p\seisuu$
for some prime $p$.
If $d$ is prime to $p$,
the derivation $t\del_t$ of $R(\!(t)\!)$
has the natural lift to $R(\!(t_d)\!)$,
which is same as $d^{-1}\cdot t_d\del_{t_d}$.
We put $I_t(g):=\sum (d/j)\cdot g_j\cdot t_d^j$
for any $g=\sum_{j\not\equiv 0\modulo p}g_j\cdot t^j
 \in R(\!(t)\!)$.
We have $t\del_t\bigl(I_t(g)\bigr)=g$
and $I_t\bigl(t\del_t g\bigr)=g$.

When $R$ is a subring of $\cnum$
finitely generated over $\seisuu$,
let $S(R,p)$ denote the set of the generic points of
the irreducible components of $\Spec(R\otimes\seisuu/p\seisuu)$
for each prime number $p$,
and we put $S(R):=\bigcup_{p}S(R,p)$.
For each $\eta\in S(R)$,
let $k(\eta)$ denote the corresponding field,
and let $\etabar\lrarr\eta$ denote a morphism
such that $k(\etabar)$ is an algebraic closure of $k(\eta)$.

We use the notation $M_r(R)$ to denote
the set of the $r$-th square matrices
over $R$, in general.

\subsection{Irregular value}

\subsubsection{Definition}

Let $k$ be a field,
whose characteristic number is denoted by $p$.
Let $E$ be a locally free $k[\![t]\!]$-module of rank $r$.
We use the notation $E(\!(t)\!)$ to denote $E\otimes k(\!(t)\!)$.
Let  $\nabla$ be a meromorphic connection of $E(\!(t)\!)$
such that $\nabla(\del_t)(E)\subset E\cdot t^{-\mu}$
for some non-negative integer $\mu$.
\begin{assumption}
\label{assumption;07.4.27.2}
If $p>0$,
we assume that $r$ and $\mu$
are sufficiently smaller than $p$,
say $10\cdot r!\cdot\mu<p$.
\hfill\qed
\end{assumption}
Let $\kuebar$ denote an algebraic closure of $k$.
Then, it is known (see \cite{andre}, for example)
that we have the unique subset
 $\Irr\bigl(E(\!(t)\!),\nabla\bigr)\subset k(\!(t_d)\!)/k[\![t_d]\!]$
and the unique decomposition
\begin{equation}
\label{eq;07.4.27.1}
  \bigl(E(\!(t)\!),\nabla\bigr)\otimes \kuebar(\!(t_d)\!)\simeq
 \bigoplus_{\gminia\in\Irr(E(\!(t)\!),\nabla)}
 \bigl(E_{\gminia}(\!(t_d)\!),
 \nabla_{\gminia}\bigr)
\end{equation}
for some appropriate factor $d$ of $r!$,
such that the following holds:
\begin{itemize}
\item
 For any element $\gminia\in\Irr(E(\!(t)\!),\nabla)$,
 take a lift $\gminiatilde\in k(\!(t_d)\!)$,
 and then $\nabla_{\gminia}-d\gminiatilde\cdot \id_{E_{\gminia}}$
 is a logarithmic connection of $E_{\gminia}$.
 The elements of $\Irr(E(\!(t)\!),\nabla)$ or their lifts
 are called the irregular values of $(E,\nabla)$.
\end{itemize}
The decomposition is called the irregular decomposition
in this paper.
We usually use the natural lifts of $\gminia$ in $k(\!(t_d)\!)_{<0}$,
and denote them by the same letter $\gminia$.
We have $\ord_t(\gminia)\geq -\mu+1$,
and $d$ is a factor of $r!$,
and hence $\ord_{t_d}(\gminia)>-p$
under the assumption \ref{assumption;07.4.27.2}.

If the irregular decomposition exists
on $k(\!(t)\!)$, 
then we say that $(E,\nabla)$ is unramified.
The following lemma easily follows from
the uniqueness of the irregular decomposition.
\begin{lem}
Let $k'$ be an algebraic extension of $k$,
and let $d'$ be a divisor of $d$.
If all the irregular values are contained in $k'(\!(t_{d'})\!)$,
then $(E,\nabla)\otimes k'(\!(t_{d'})\!)$ is unramified.
\hfill\qed
\end{lem}

\subsubsection{Connection form of Deligne-Malgrange lattice}

We have another characterization of
the irregular values.
For simplicity, we assume that 
$\bigl(E(\!(t)\!),\nabla\bigr)$ is unramified
and that $k$ is algebraically closed.

\begin{df}
We say that $E$ is a Deligne-Malgrange lattice of $E(\!(t)\!)$,
if the irregular decomposition {\rm(\ref{eq;07.4.27.1})}
is given on $k[\![t]\!]$ not only on $k(\!(t)\!)$,
i.e.,
$E=\bigoplus E_{\gminia}$.

If $E$ is Deligne-Malgrange,
we have the logarithmic connection
 $\nabla^{\reg}=\bigoplus
 \nabla^{\reg}_{\gminia}$,
where $\nabla^{\reg}_{\gminia}:=
 \nabla_{\gminia}-d\gminia\cdot\id_{E_{\gminia}}$.
We say $E$ is a strict Deligne-Malgrange lattice,
if $\alpha-\beta$ are not integers
for any two distinct eigenvalues $\alpha,\beta$ of
$\Res(\nabla^{\reg})$.
\hfill\qed
\end{df}

Let $\vecv$ be any frame of $E$.
Let $A\in M_r\bigl(k(\!(t)\!)\bigr)$ be determined by
$\nabla(t\del_t)\vecv=\vecv\cdot A$.
Let $\Sp(A)\in k(\!(t_d)\!)$ denote the set of
the eigenvalues of $A$ for some appropriate $d$.
For any $\alpha\in\Sp(A)$,
we have the negative part 
$\alpha_-\in k(\!(t_d)\!)_{<0}$
and $I_t(\alpha_-)\in k(\!(t_d)\!)_{<0}$
as explained in Subsection \ref{subsection;07.4.26.1}.

\begin{lem}
\label{lem;07.4.26.3}
If $E$ is Deligne-Malgrange,
we have
$\Irr\bigl(E(\!(t)\!),\nabla\bigr)=
 \bigl\{I_t(\alpha_-)\,\big|\,
 \alpha\in \Sp(A)\bigr\}$.
\end{lem}
\pf
We take a frame $\vecv_1$ of $E$
compatible with the irregular decomposition,
and $A_1$ is determined as above.
Then, we have $A_1$ has the decomposition corresponding
to the irregular decomposition,
$ A_1=
\bigoplus  
 (t\del_t\gminia+R_{\gminia})$,
where $R_{\gminia}\in M_r(k[\![t]\!])$.
Hence
the claim of the lemma clearly holds
for the frame $\vecv_1$.

For any frame $\vecv$ of $E$,
we have $G\in \GL\bigl(k[\![t]\!]\bigr)$
such that $\vecv=\vecv_1\cdot G$.
We have the relation
$A=G^{-1}\cdot A_1\cdot G+G^{-1}\cdot t\del_t G$,
i.e.,
$G\cdot A\cdot G^{-1}
=A_1+(t\del_tG)\cdot G^{-1}$,
where $t\del_t G\cdot G^{-1}\in M_r(k[\![t]\!])$.

Hence, the claim is reduced to the following general lemma.
\begin{lem}
\label{lem;07.4.26.2}
Let $\Gamma\in M_r(k[\![t]\!])$ be a diagonal matrix
whose $(i,i)$-entry is given by $\alpha_i$.
Let $B$ be any element of $t^m\cdot M_r(k[\![t]\!])$
for a positive integer $m>0$.
Then, any eigenvalue $\beta\in k[\![t_d]\!]$
of $\Gamma+B$
satisfies $\ord_t(\beta-\alpha_i)\geq m$
for some $\alpha_i$.
\end{lem}
\pf
Let $e_1,\ldots,e_r$ denote the canonical base
of $k(\!(t)\!)^r$.
Let $v=\sum f_i\cdot e_i$ be an eigenvector
of $\Gamma+B$
corresponding to the eigenvalue $\beta$.
We may assume $\ord_t(f_{i_0})=0$ for some $i_0$.
We obtain
$\ord_t\bigl( (\alpha_i-\beta)\cdot f_i\bigr)\geq m$ for any $i$,
and hence $\ord_t(\alpha_{i_0}-\beta)\geq m$.
Thus,
we obtain Lemma \ref{lem;07.4.26.2}
and Lemma \ref{lem;07.4.26.3}.
\hfill\qed

\subsubsection{$p$-curvature}

In the case $p>0$,
we have the other characterization of the irregular values.
For simplicity, we assume $k=\kuebar$.
Let $\Fr:k(\!(t)\!)\lrarr k(\!(t)\!)$ be
the absolute Frobenius morphism,
i.e.,
$\Fr(f)=f^p$.
Applying $\Fr$ to the coefficients,
we obtain the homomorphism
$k(\!(t)\!)[T]\lrarr k(\!(t)\!)[T]$,
which is also denoted by $\Fr$.
Let $\psi$ be the $p$-curvature of $\nabla$.
(See \cite{katz1} and \cite{katz2}, for example).
Due to the observation of Bost-Laszlo-Pauly
(\cite{laszlo-pauly}),
there exists a polynomial $P_{\nabla}(T)\in k(\!(t)\!)[T]$
of degree $r$,
such that
$\det\bigl(T-\psi(t\del_t)\bigr)
=\Fr\bigl(P_{\nabla}\bigr)(T)$.
Let $\Sol(P_{\nabla})$ denote the set of
the solutions of $P_{\nabla}(T)=0$.
Then $\Sol(P_{\nabla})\subset
 k(\!(t_d)\!)$ for some appropriate factor $d$ of $r!$.
Because of $\nabla(\del_t)(E)\subset E\cdot t^{-\mu}$,
we have $\psi(\del_t)(E)\subset E\cdot t^{-\mu\cdot p}$.
Hence we have $\ord_t(\alpha)\geq -\mu+1$
for any solution $\alpha\in \Sol(P_{\nabla})$.
Under the assumption \ref{assumption;07.4.27.2},
we obtain $\ord_{t_d}(\alpha)>-p$
for any $\alpha\in\Sol(P_{\nabla})$.

\begin{lem}
\label{lem;07.4.14.3}
Under the assumption
{\rm\ref{assumption;07.4.27.2}},
we have $\Irr(E(\!(t)\!),\nabla)
=\bigl\{I_t(\alpha_-)\,\big|\,
 \alpha\in \Sol(P_{\nabla})\bigr\}$.
\end{lem}
\pf
We may assume that $(E,\nabla)$ is unramified
and Deligne-Malgrange.
Hence, we have only to consider the case 
where $(E,\nabla)$ has the unique irregular value,
i.e.,
$\nabla=d\gminia\cdot \id_{E}+\nabla^{\reg}$,
where $\gminia\in k(\!(t)\!)_{<0}$,
$\ord_{t}(\gminia)>-p$,
and $\nabla^{\reg}$ is logarithmic.
Let $\psi_{\reg}$ denote the $p$-curvature of $\nabla^{\reg}$.
By a general formula
(\cite{robba-christol}. See also 
 Lemma 3.4 of \cite{tsuchimoto}),
we have
$\psi(t\del_t)=\psi_{\reg}(t\del_t)+(t\del_t\gminia)^p$,
where $\psi_{\reg}(t\del_t)\in M_r(k[\![t]\!])$.
Then the claim of the lemma follows from
Lemma \ref{lem;07.4.26.2}.
\hfill\qed

\subsection{Preliminary from elementary algebra}

The following arguments are standard and well known.
We would like to be careful about some finiteness,
and we give just an outline.
Let $R$ be a regular ring.
Let $P_t(T)\in R[\![t]\!][T]$ be a monic polynomial:
$P_t(T)=T^r+\sum_{j=0}^{r-1}a_j(t)\cdot T^j$.
The specialization at $t=0$ is denoted by $P_0(T)$.
\begin{lem}
\label{lem;07.4.14.1}
Assume that $P_0(T)=\huebar_1(T)\cdot \huebar_2(T)$ in $R[T]$
such that $\huebar_1$ and $\huebar_2$ are monic polynomials
and coprime in $K[T]$.
Then, we have the decomposition
$P(T)=h_1(T)\cdot h_2(T)$ in $R'[\![t]\!][T]$,
where $R'$ is the localization of $R$
with respect to some $f_1,\ldots,f_m\in R$
depending on $P_0(T)$,
and $h_i(T)$ are monics such that
$h_i(T)_{|t=0}=\huebar_i(T)$.
\end{lem}
\pf
There exist $F_i\in K[T]$ $(i=1,2)$
such that $1=\huebar_1(T)\cdot F_1(T)+\huebar_2(T)\cdot F_2(T)$.
We may take a finite localization
$R'$ of $R$ so that $F_i\in R'[T]$.
For any $Q(T)\in R'[T]$ such that $\deg_T Q<r$,
we have
$\huebar_1\cdot (F_1Q)+\huebar_2\cdot (F_2Q)=Q$.
Take $H,G\in R'[T]$ such that
$\deg_T(H)<\deg_T(\huebar_1)$ and
$F_2Q=\huebar_1\cdot G+H$.
We put $\alpha=F_1Q+G\huebar_2$,
and then we have
$\huebar_1\cdot \alpha+\huebar_2\cdot H=Q$.
Note
$\deg_T(\huebar_1)+\deg_T(\huebar_2)=\deg_T P_0=r$.
Then, we have
$\deg_T(\alpha)+\deg_T \huebar_1
\leq
 \max\bigl(\deg_T Q,\deg_T \huebar_2+\deg_T H\bigr)
<r$.
Hence,
$\deg_T(\alpha)<r-\deg_T(\huebar_1)=\deg_T(\huebar_2)$.

Assume we are given
$h_{a,j}(T)$ $(a=1,2,\,\,j=1,\ldots,L)$
such that $\deg_T h_{1,j}<\deg_T (\huebar_1)$,
$\deg_T h_{2,j}< \deg_T(\huebar_2)$
and 
\[
 \left(
 \huebar_1(T)+\sum_{j=1}^Lh_{1,j}(T)t^j
 \right)
\cdot
 \left(
 \huebar_2(T)+\sum_{j=1}^Lh_{2,j}(T)t^j
 \right)
-\sum_{j=0}^LP_{j}(T)t^j
\equiv 0\mod t^{L+1} 
\]
By using the above remark,
it is easy to show that we can take
$h_{a,L+1}$ $(a=1,2)$ such that
$\deg_T h_{1,L+1}<\deg_T(\huebar_1)$,
$\deg_T h_{2,L+1}< \deg_T(\huebar_2)$
and 
\[
 \left(
 \huebar_1(T)+\sum_{j=1}^{L+1}h_{1,j}(T)t^j
 \right)
\cdot
 \left(
 \huebar_2(T)+\sum_{j=1}^{L+1}h_{2,j}(T)t^j
 \right)
-\sum_{j=0}^{L+1}P_{j}(T)t^j
\equiv 0\mod t^{L+2} 
\]
Thus, by an inductive argument,
we can construct the desired $h_1$ and $h_2$.
\hfill\qed

\begin{lem}
\label{lem;07.4.27.5}
Let $P_t(T)\in R[\![t]\!][T]$ (resp. $R(\!(t)\!)[T]$) be a monic polynomial.
There exists an appropriate number $e$,
such that it has a roots in $R'[\![t_e]\!]$ (resp. $R'(\!(t_e)\!)$)
where $R'$ is obtained from $R$
by finite algebraic extensions and localizations.
\end{lem}
\pf
Let $P_t(T)=\sum_{j=0}^n a_j(t)\cdot T^j$.
We may assume that $n!$ is invertible in $R$.
Let $\nu(P_t)$ denote the number
$\min_j\bigl\{
 \ord_t(a_j)\big/(n-j)\bigr\}$.
We use the induction on the numbers $\deg_TP_t$
and $\nu(P_t)$.
For simplicity, we use $\nu$ instead of $\nu(P_t)$,
and let $d$ be the minimal positive integer
such that $\nu\in d^{-1}\cdot\seisuu$.
We formally use the notation 
$t^{\nu}$ to denote $t_{d}^{d\cdot\nu}$.
We have the following monic polynomial:
\[
 Q_t(T'):=t^{-n\nu}P_t(t^{\nu}T')
=\sum_{j=0}^n a_j(t)t^{-(n-j)\nu}T^{\prime\,j}
=\sum_{j=0}^nb_j(t)T^{\prime\,j}
\in R[\![t_{d}]\!][T']
\]
We have
$ d^{-1}\cdot\ord_{t_d}(b_j)
=\ord_t(a_j)-(n-j)\cdot\nu\geq 0$,
and we have $\ord_{t_{d}}(b_{j_0})=0$ for some $j_0$.
We put $Q_0(T')=\sum_{j=0}^nb_j(0)T^{\prime\,j}
 \in R[T']$.

\vspace{.1in}
\noindent
{\bf Case 1}\,\,
Assume $Q_0(T')$ has at least two different roots.
Then, there exists a finite algebraic extension $R_1$ of $R$
such that we have the decomposition
$Q_0(T')=\huebar_1(T')\huebar_2(T')$ in $R_1[T']$,
and $\huebar_1$ and $\huebar_2$ are coprime
in $K_1[T']$,
where $K_1$ denotes the quotient field of $R_1$.
Because of Lemma \ref{lem;07.4.14.1},
we have
$Q_t(T')=h_1(T')\cdot h_2(T')$
in $R_1'[\![t_d]\!][T']$,
where $R_1'$ is a localization of $R_1$
with respect to some finite elements.
By the hypothesis of the induction on the degree
with respect to $T'$,
$h_i(T')$ $(i=1,2)$ have the roots $\alpha$ in $R_2[\![t_{d}]\!]$,
where $R_2$ is obtained from  $R_1'$
by finite algebraic extensions and localizations.
And $t^{\nu}\alpha$ give the roots of $P_t(T)$.

\vspace{.1in}
\noindent
{\bf Case 2}\,\,
In the case $Q_0(T')=(T'-\alpha)^n$,
we have $n\alpha\in R$,
and hence $\alpha\in R$.
We have
 $\ord_t(a_n)/n=\ord_t(a_{n-1})/(n-1)=\nu$,
and hence $\nu\in\seisuu$ and $d=1$.
We put 
$H_t(T):=P_t(T+t^{\nu}\alpha)
=\sum_{j=0}^n c_j(t)T^j$.
We have $\min\bigl(\ord(c_j)(n-j)^{-1}\bigr)>\nu$.

We continue the process.
If we reach the case 1,
we can reduce the degree
with respect to $T$.
If we do not reach the case $1$,
it is shown that
$P_t(T)=(T-\gminia)^n$ for some $\gminia\in R[\![t]\!]$
(resp. $\gminia\in R(\!(t)\!)$).
Thus we are done.
\hfill\qed

\begin{cor}
Any $P(s,t)(T)\in R[\![s,t]\!][T]$
has the roots in $R'_P(\!(s_d)\!)[\![t_d]\!]$,
and any $P(s,t)(T)\in R[\![s]\!](\!(t)\!)[T]$
has the roots in $R'_P(\!(s_d)\!)(\!(t_d)\!)$.
Here $R'_P$ is obtained from $R$,
depending on $P$,
by finite algebraic extensions and localizations,
and $d$ denotes an appropriate positive integer.
\hfill\qed
\end{cor}

Let $R$ be an integral domain such that
$\seisuu\subset R$.
Let $K$ denote the quotient field of $R$.
Let $(\nbige,\nabla)$ be a meromorphic connection
on $R(\!(t)\!)$.
\begin{lem}
There exists an extension $R'$
obtained from $R$ by
finite algebraic extensions and localizations,
with the following property:
\begin{itemize}
\item
 The irregular values of $(\nbige,\nabla)\otimes K(\!(t)\!)$
 are contained in $R'(\!(t_d)\!)$.
\item
 The irregular decomposition
 and a Deligne-Malgrange lattice 
 are defined on $R'(\!(t_d)\!)$.
\end{itemize}
\end{lem}
\pf
We need only a minor modification for the argument
given in \cite{levelt},
and hence we give just an outline.
We put $D=\nabla(t\del_t)$.
Let $\nbigk$ be the quotient field of $R(\!(t)\!)$.
By applying the argument of Deligne \cite{d}
to $E\otimes\nbigk$ with the derivation $D$,
we can take $e\in \nbige$
such that $e,D(e),\ldots, D^{r-1}(e)$
give a base of the $\nbigk$-vector space $E\otimes\nbigk$.
We have the relation
$D^re+\sum_{j=0}^{r-1}a_j\cdot D^je=0$
where $a_j\in \nbigk$.
There exists a finite localization $R_1$ of $R$
such that $a_j\in R_1(\!(t)\!)$.

We put $\nu:=\min_j\bigl\{\ord_t(a_j)\big/(r-j)\bigr\}$.
Note that $\nu\geq 0$ implies the regularity of 
the connection.
Let $d$ denote the minimal positive integer
such that $d\cdot\nu\in \seisuu$.
We put $f_{i+1}:=t^{-\nu\cdot i} D^ie$ $(i=0,\ldots,r-1)$,
and $\vecf=(f_i\,|\,i=1,\ldots,r)$.
Let $A\in M_r\bigl(R_1(\!(t_d)\!)\bigr)$ be determined by
$D\vecf=\vecf A$.
Then, $A$ is of the form $t^{\nu}\bigl(A_0+t_d\cdot A_1(t_d)\bigr)$
such that 
(i) $A_1\in M_r\bigl(R_1[\![t_d]\!]\bigr)$,
(ii) $A_0\in M_r(R_1)$ whose $(i,j)$-entries are as follows:
\[
 \bigl(A_0\bigr)_{i,j}=
 \left\{ \begin{array}{ll}
 1 & (i=j+1)\\
 -(t^{(-r+i-1)\nu}a_{i-1})_{|t_d=0} & (j=r)\\
 0 &(\mbox{\rm otherwise})
 \end{array}
 \right.
\]
By the choice, one of $(A_0)_{i,r}$ not $0$.

\noindent
{\bf Case 1}\,
Let us consider the case where $A_0$ has
at least two distinct eigenvalues.
There exists a finite extension $R_2$ such that 
(i) we have $G\in \GL_r(R_2)$ for which $G^{-1}A_0G$ is Jordan,
(ii) the difference of any two distinct eigenvalues of 
 $A_0$ are invertible in $R_2$.
By a standard argument (See \cite{levelt}, or 
 the proof of Lemma \ref{lem;07.9.19.13} below),
we can show that
there exists $G_1\in \GL_r\bigl(R_2[\![t_d]\!]\bigr)$
such that
(i) $G_{1\,|\,t_d=0}=G$,
(ii) let $\vecg=\vecf\cdot G_1$
 and $D\vecg=\vecg\cdot B$,
 then $B$ is decomposed into a direct sum
 of matrices with smaller sizes.
Hence, we obtain a decomposition
into connections with lower ranks.
Thus, we can reduce the problem
to the lower rank case.

\vspace{.1in}
\noindent
{\bf Case 2}\,
If $A_0$ has the unique eigenvalues $\alpha\in R_1$,
it can be shown that $d=1$ and $\nu\in\seisuu_{<0}$,
as in the proof of Lemma \ref{lem;07.4.27.5}.
We put $\nabla'=\nabla-t^{\nu}\alpha\cdot dt/t$
and $D'=\nabla'(t\del_t)$.
Let $\nbigk_1$ be the quotient field of $R_1(\!(t)\!)$.
It can be shown that
$e,D'e,\ldots, (D')^{r-1}e$ give a base of
$\nbige\otimes \nbigk_1$.
Let $a'_j$ be determined by 
$D^{\prime\,r}e+
 \sum a_j'\cdot D^{\prime\,j}e=0$.
Then, we have $a'_j\in R_1(\!(t)\!)$
and $\nu(\nabla')=\min\bigl\{ \ord_t(a'_j)\big/(r-j) \bigr\} \geq \nu+|\nu|/r$.
We continue the process.
After the finite steps,
we will arrive at the case $1$
or the case $\nu(\nabla')\geq 0$.
\hfill\qed

\vspace{.1in}

\begin{cor}
\label{cor;07.4.28.1}
Let $(\nbige,\nabla)$ be a meromorphic connection
on $R[\![s]\!](\!(t)\!)$.
Then, there exists an extension $R'$,
which is obtained from $R$ by
finite algebraic extensions and localizations,
with the following property:
\begin{itemize}
\item
 The irregular values of $(\nbige,\nabla)\otimes K(\!(s)\!)(\!(t)\!)$
 are contained in $R'(\!(s_d)\!)(\!(t_d)\!)$.
\item
 The irregular decomposition 
 and a Deligne-Malgrange lattice 
 are defined on $R'(\!(s_d)\!)(\!(t_d)\!)$.
\hfill\qed
\end{itemize}
\end{cor}

\subsection{Good formal structure}
\label{subsection;07.4.30.1}

Let $X$ be a complex algebraic surface,
with a simple normal crossing divisor $D$.
Let $(\nbige,\nabla)$ be a meromorphic flat connection
on $(X,D)$.
We recall the notion of good formal structure,
by following \cite{sabbah4}.

If $P$ is a smooth point of $D$,
we take a holomorphic coordinate $(U,t,s)$ around $P$
such that $t^{-1}(0)=U\cap D$.
For a positive integer $d$,
we take a ramified covering $\varphi_d:U_d\lrarr U$
given by $(t_d,s)\longmapsto (t_d^d,s)$.
We put $D_d:=\{t_d=0\}\subset U_d$.
Let $M(U_d,D_d)$ (resp. $H(U_d)$)
denote the space of meromorphic 
(resp. holomorphic) functions on $U_d$
whose poles are contained in $D_d$.
For any element $\gminia$ of $M(U_d,D_d)/H(U_d)$,
we have the natural lift to $M(U_d,D_d)$
which is also denoted by $\gminia$.
Let $\Dhat_d$ denote the formal space
obtained as the completion of $U_d$ along $D_d$.
(See \cite{bingener}, for example.)

\begin{df}
We say that $(\nbige,\nabla)$ has the good formal structure at $P$,
if the following holds for some $(U,t,s)$ and
some $d\in\seisuu_{>0}$:
\begin{itemize}
\item
We have the finite subset
$\Irr(\nbige,\nabla)\subset M(U_d)/H(U_d)$
and the decomposition:
\[
 \varphi_d^{\ast}(\nbige,\nabla)_{|\Dhat_d}
=
\bigoplus_{\gminia\in \Irr(\nbige,\nabla)}
 \bigl(\nbige_{\gminia},\nabla_{\gminia}\bigr)
\]
Here
$\nabla_{\gminia}^{\reg}:=
 \nabla_{\gminia}-d\gminia\cdot \id_{\nbige_{\gminia}}$
are regular.
\item
 For any non-zero $\gminia\in\Irr(\nbige,\nabla)$,
 the $0$-divisor of $\gminia$  has no intersection
 with $D_d$.
\item
 For any two distinct $\gminia,\gminib\in\Irr(\nbige,\nabla)$,
 the $0$-divisor of $\gminia-\gminib$
 has no intersection with $D_d$.
\hfill\qed
\end{itemize}
\end{df}

If $P$ is a cross point of $D$,
we take a holomorphic coordinate $(U,t,s)$
such that $D\cap U=\{t\cdot s=0\}$.
For each $d\in\seisuu_{>0}$,
we take a ramified covering
$\varphi_d:U_d\lrarr U$ given by
$(t_d,s_d)\longmapsto (t_d^d,s_d^d)$.
We put $D_d:=\{t_d\cdot s_d=0\}$ and $P_d:=(0,0)$.
Let $\Phat_d$ denote the formal space
obtained as the completion of $U_d$ at $P_d$.

Let $M(U_d,D_d)$ (resp. $H(U_d)$) denote 
the space of the meromorphic (holomorphic) functions
on $U_d$ whose poles are contained in $D_d$.
For any element $\gminia$ of $M(U_d,D_d)/H(U_d)$,
we have the natural lift to $M(U_d,D_d)$,
which is also denoted by $\gminia$.

We use the partial order $\leq_{\seisuu^2}$ on $\seisuu^2$
given by $(a_1,a_2)\leq_{\seisuu^2} (a_1',a_2')
 \Longleftrightarrow a_i\leq a_i'\,\,(i=1,2)$.
For any element $f=\sum f_{i,j}\cdot s^i\cdot t^j\in M(U_d,D_d)$,
let $\ord(f)$ denote the minimum of the set
$\min\bigl\{(i,j)\,\big|\,f_{i,j}\neq 0\bigr\}$,
if it exists.

\begin{df}
We say that $(\nbige,\nabla)$ has the good formal structure
if the following holds:
\begin{itemize}
\item
 We have the finite subset $\Irr(\nbige,\nabla)
 \subset M(U_d)/H(U_d)$
 and the decomposition for some $d\in\seisuu_{>0}$:
\[
 \varphi_d^{\ast} (\nbige,\nabla)_{|\Phat_d}
=
 \bigoplus_{\gminia\in \Irr(\nbige,\nabla)}
 \bigl(\nbige_{\gminia},\nabla_{\gminia}\bigr)
\]
Here
$\nabla_{\gminia}^{\reg}:=
 \nabla_{\gminia}-d\gminia\cdot \id_{\nbige_{\gminia}}$
are regular.
\item
$\ord(\gminia)$ exists in $\seisuu_{\leq 0}^2-\{(0,0)\}$
for each non-zero $\gminia\in\Irr(\nbige,\nabla)$.
\item
$\ord(\gminia-\gminib)$ exists in $\seisuu_{\leq 0}^2-\{(0,0)\}$
for any two distinct $\gminia,\gminib\in\Irr(\nbige,\nabla)$.
And the set
$\bigl\{\ord(\gminia-\gminib)\,|\,
 \gminia,\gminib\in\Irr(\nbige,\nabla)\bigr\}$
is totally ordered with respect to
 the above order $\leq_{\seisuu^2}$.
\hfill\qed
\end{itemize}
\end{df}

\begin{df}
A point $P$ is called turning with respect to $(\nbige,\nabla)$,
if $(\nbige,\nabla)$ does not have a good formal structure
at $P$.
\hfill\qed
\end{df}

\subsection{A sufficient condition
for the existence of the good formal structure}

\subsubsection{Preliminary}

Let $E$ be a free $\cnum[\![s,t]\!]$-module.
Let $\nabla_t:E\lrarr E\otimes
 \Omega^1_{
   \cnum[\![s]\!](\!(t)\!)/\cnum[\![s]\!]}(\ast st)$
be a connection such that
the following holds for some $k\geq 1$ and $p\geq 0$:
\[
 \nabla_t(t^{k+1}s^p\del_t)(E)\subset E
\]
In that case,
$\nabla_t\bigl(t^{k+1}s^p\del_t\bigr)$
induces the endomorphism of
$E_0:=E_{|t=0}$,
which is denoted by $F_0$.

\begin{lem}
\label{lem;07.9.19.12}
If $F_0$ is invertible,
any meromorphic flat section 
$f=\sum _{j\geq -N}f_j\cdot t^j$ of $E$
is $0$.
\end{lem}
\pf
Let $f$ be a meromorphic flat section of $E$.
Assume $f\neq 0$.
We may assume that
$-N=\min\{j\,|\,f_j\neq 0\}$.
From $\nabla(t^{k+1}s^p\del_t)f=0$,
we have $F_0(f_{-N})=0$.
Because $F_0$ is invertible,
we obtain $f_{-N}=0$,
which contradicts with the choice of $N$.
\hfill\qed

\begin{lem}
\label{lem;07.9.19.13}
Assume the following:
\begin{itemize}
\item
 We have the decomposition
 $(E_0,F_0)=
 (E_0^{(1)},F_0^{(1)})
 \oplus (E_0^{(2)},F_0^{(2)})$.
\item
 The eigenvalues of $F_0^{(i)}$ are
 contained in $\cnum[\![s]\!]$.
 If $\gminib_i$ $(i=1,2)$ are the eigenvalues of
 $F_0^{(i)}$,
 we have $(\gminib_1-\gminib_2)_{|s=0}\neq 0$.
\end{itemize}
Then, we have the unique $\nabla_t$-flat decomposition
$E=E^{(1)}\oplus E^{(2)}$ such that
the restriction to $t=0$ is the same as 
$E_0=E^{(1)}_0\oplus E^{(2)}_0$.
\end{lem}
\pf
We closely follow the argument in \cite{levelt}.
Let $\vecv$ be a frame of $E$
such that $\vecv_{|t=0}$ is compatible with
the decomposition $E_0=E^{(1)}_0\oplus E^{(2)}_0$.
Then, $\vecv$ is divided as $(\vecv^{(1)},\vecv^{(2)})$,
where $\vecv^{(i)}_{|t=0}$ are the frames of
$E^{(i)}_0$.
Let $A=\sum_{j=0}^{\infty} A_j(s)\cdot t^j$ be 
determined by the following:
\[
 \nabla(t^{k+1}s^p\del_t)\vecv
=\vecv\cdot A.
\]
We have the following decomposition
corresponding to the decomposition of the frame
$\vecv=(\vecv^{(1)},\vecv^{(2)})$:
\begin{equation}
 A_j=\left(\begin{array}{cc}
 A_j^{(11)} & A_j^{(12)}\\
 A_j^{(21)} & A_j^{(22)}
 \end{array}
 \right)
\end{equation}
By the assumption,
we have $A_0^{(12)}=0$ and $A_0^{(21)}=0$.
For the change of the frame 
from $\vecv$ to $\vecv\cdot G$,
we have the following:
\[
 \nabla(t^{k+1}s^p\del_t)\bigl(\vecv\cdot G\bigr)
=\bigl(\vecv\cdot G\bigr)\cdot \Atilde(G),
\quad
\Atilde(G):=
G^{-1}AG+t^{k+1}s^pG^{-1}\del_tG
\]
We consider the formal transform $G$
of the following form:
\[
 G=I+\left(
 \begin{array}{cc}
 0 & X\\
 Y & 0
 \end{array}
 \right),
\quad
 X=\sum_{j=1}^{\infty}X_j(s)\cdot t^j,
\quad
 Y=\sum_{j=1}^{\infty}Y_j(s)\cdot t^j
\]
Here the entries of $X_j(s)$ and $Y_j(s)$
are contained in $\cnum[\![s]\!]$.
We want to determine $X_j$ and $Y_j$
by the following condition:
\begin{itemize}
\item
 The $(1,2)$-component and the $(2,1)$-component
 of $\Atilde(G)$
 are $0$.
\item
 The $(1,1)$-component of $\Atilde(G)$
 is of the form
 $A^{(11)}_0+B^{(11)}$,
 where the entries of
 $B^{(11)}$ are contained in $t\cdot\cnum[\![s,t]\!]$.
 Similarly,
 The $(2,2)$-component is of the form
 $A^{(22)}_0+B^{(22)}$,
 where the entries of
 $B^{(22)}$ are contained in $t\cdot\cnum[\![s,t]\!]$.
\end{itemize}
We obtain the following equations
for $Y$ and $B^{(11)}$:
\[
 A^{(11)}+A^{(12)}Y-A_0^{(11)}-B^{(11)}=0,
\quad
 A^{(21)}+A^{(22)}Y+t^{k+1}s^p\del_tY
-Y(A_0^{(11)}+B^{(11)})=0
\]
Then, we obtain the following equation for $Y$:
\[
 A_0^{(22)}Y-YA_0^{(11)}
-Y(A^{(11)}-A_0^{(11)})
+(A^{(22)}-A_0^{22})Y
-YA^{(12)}Y
+t^{k+1}s^p\del_tY
+A^{(21)}=0
\]
For the expansion
$Y=\sum_{j=1}^{\infty} Y_j(s)\cdot t^j$,
we obtain the following equations:
\begin{equation}
 \label{eq;07.9.19.11}
 A_0^{(22)}Y_j-Y_jA_0^{(11)}
-\sum_{\substack{l+m=j\\l,m\geq 1}}
 Y_lA^{(11)}_m
+\sum_{\substack{l+m=j\\l,m\geq 1}}
 A_l^{(22)}Y_m
-\sum_{\substack{l+m+n=j,\\l,m,n\geq 1}}
 Y_lA_m^{(12)}Y_n
+(j-k)s^pY_{j-k}\cdot \chi_{j\geq k}
+A_j^{(21)}=0
\end{equation}
Here $\chi_{j\geq k}=0$ if $j<k$
and $\chi_{j\geq k}=1$ if $j\geq k$.
When we are given $Y_m$ $(1\leq m\leq j-1)$
whose entries are contained in $\cnum[\![s]\!]$,
we have the unique solution $Y_j$
of (\ref{eq;07.9.19.11}),
whose entries are contained in
$\cnum[\![s]\!]$.
Hence, we have appropriate $Y$ and $B^{(11)}$.
Similarly, we have appropriate $X$ and $B^{(22)}$.
Thus, we can conclude the existence of
the desired decomposition $E=E^{(1)}\oplus E^{(2)}$.
The uniqueness follows from
Lemma \ref{lem;07.9.19.12}.
\hfill\qed

\vspace{.1in}
Let us consider the case where
$\nabla_t$ comes from 
a flat meromorphic connection 
$\nabla:E\lrarr
 E\otimes\Omega^{1}_{
   \cnum[\![s]\!](\!(t)\!)/\cnum}(\ast st)$.

\begin{lem}
\label{lem;07.9.19.14}
Assume the hypothesis in Lemma {\rm\ref{lem;07.9.19.13}}.
The decomposition $E=E^{(1)}\oplus E^{(2)}$
is $\nabla$-flat.
\end{lem}
\pf
We may assume $\vecv=(\vecv_1,\vecv_2)$ is compatible with
the decomposition $E=E^{(1)}\oplus E^{(2)}$.
Let $A$ and $B$ be determined by the following:
\[
 \nabla(t^{k+1}\del_t)\vecv=\vecv\cdot A,
\quad
 A=\left(
 \begin{array}{cc}
 A^{(11)} & 0\\
 0 & A^{(22)}
 \end{array}
 \right)
\]
\[
 \nabla(\del_s)\vecv=\vecv\cdot B,
\quad
 B=\left(
 \begin{array}{cc}
 B^{(11)} & B^{(12)}\\ B^{(21)} & B^{(22)}
 \end{array}
 \right)
\]
From the relation
$\bigl[\nabla(\del_s),\nabla(t^{k+1}\del_t)\bigr]$,
we have the following equation for $B^{(12)}$:
\[
 A^{(11)}B^{(12)}
-B^{(12)}A^{(22)}
+t^{k+1}\del_tB^{(12)}=0
\]
Assume $B^{(12)}\neq 0$.
We have the expression
$B^{(12)}=\sum_{j\geq -N}B_j^{(12)}\cdot t^j$,
and we may assume $B_{-N}^{(12)}\neq 0$.
But, we have the relation
$B_{-N}^{(12)}A_0^{(11)}
-A_0^{(22)}B_{-N}^{(12)}=0$,
and hence $B_{-N}^{(12)}=0$.
Thus, we arrive at the contradiction,
and we can conclude $B^{(12)}=0$.
Similarly,
we obtain $B^{(21)}=0$.
\hfill\qed

\subsubsection{A condition}

Let $E$ be a free $\cnum[\![s,t]\!]$-module
with a flat meromorphic connection
$\nabla:E\lrarr
 E\otimes\Omega^{1}_{
   \cnum[\![s]\!](\!(t)\!)/\cnum}(\ast t)$.
We have the induced relative connection
$\nabla_t:E\lrarr
 E\otimes\Omega^{1}_{
   \cnum[\![s]\!](\!(t)\!)/\cnum[\![s]\!]}(\ast t)$.

We put $\gbigk:=\cnum(\!(s)\!)(\!(t)\!)$.
We put $(\nbige_{\gbigk},\nabla_{\gbigk}):=
 (E,\nabla_t)\otimes\gbigk$
and $E_{\gbigk}:=E\otimes\cnum(\!(s)\!)[\![t]\!]$.
We assume that
$E_{\gbigk}$ is a strict Deligne-Malgrange lattice.
The intersection of
$E_{\gbigk}$ and $E(\ast t)$
in $\nbige_{\gbigk}$ is the same as $E$,
which gives a characterization of $E$.

\begin{prop}
\label{prop;07.9.20.1}
Assume the following:
\begin{itemize}
\item
$\Irr(\nbige_{\gbigk},\nabla_{\gbigk})$
is contained in
$\cnum[\![s]\!](\!(t)\!)/\cnum[\![s,t]\!]$.
\item
$\Irr(\nbige_{\gbigk},\nabla_{\gbigk})$ is good,
in the following sense:
 \begin{itemize}
 \item Let $\gminia\in \Irr(\nbige_{\gbigk},\nabla_{\gbigk})$.
 For the expression 
 $\gminia=\sum_{j\geq \ord_t(\gminia)}\gminia_j(s)\cdot t^j$,
 we have $\gminia_{\ord_t(\gminia)}(0)\neq 0$.
 \item
 Similarly, we have
 $(\gminia-\gminib)_{\ord_t(\gminia-\gminib)}(0)\neq 0$
 for any two distinct 
 $\gminia,\gminib\in \Irr(\nbige_{\gbigk},\nabla_{\gbigk})$.
 \end{itemize}
\end{itemize}
Then, $\bigl(E(\ast t),\nabla\bigr)$ has the good formal structure.
\end{prop}
\pf
We put $k(E):=-\min\bigl\{\ord_t(\gminia)\,\big|\,
 \gminia\in\Irr(\nabla_{\gbigk})\bigr\}$.
Assume $k(E)\geq 1$.
Because $E_{\gbigk}$ is a Deligne-Malgrange
lattice of $\nbige_{\gbigk}$,
we have $\nabla(t^{k(E)+1}\del_t)E\subset E$.
Let $F_0$ denote the endomorphism
of $E_{|0}$ induced by $\nabla(t^{k(E)+1}\del_t)$.
The eigenvalues of $F_0$ are given by
$(t^{k(E)+1}\del_t\gminia)_{t=0}$
($\gminia\in \Irr(\nabla_{\gbigk})$).
By using Lemma \ref{lem;07.9.19.13}
and Lemma \ref{lem;07.9.19.14},
we obtain the decomposition:
\[
 (E,\nabla)
=\bigoplus_{\gminib\in S}
 (E_{\gminib},\nabla_{\gminib}),
\quad
 S:=\Bigl\{
 \gminib= 
 \bigl(
t^{k(E)}\gminia
\bigr)_{|t=0}\in\cnum[\![s]\!]
 \,\Big|\,
 \gminia\in\Irr(\nabla_{\gbigk})
 \Bigr\}
\]
We put
$\nabla'_{\gminib}:=\nabla_{\gminib}
-d\bigl(t^{-k(E)}\gminib\bigr)$.
Then, $(E_{\gminib},\nabla'_{\gminib})$
also satisfy the assumption of this lemma,
and we have $k(E_{\gminib})\leq k(E)-1$.
If $k(E_{\gminib})\geq 1$,
we may apply the above argument
to $(E_{\gminib},\nabla'_{\gminib})$.
By the inductive argument,
we obtain the flat decomposition
$ (E,\nabla)
=\bigoplus_{\gminia\in \Irr(\nabla_{\gbigk})}
 (E_{\gminia},\nabla_{\gminia})$
such that 
$\nabla^{\reg}_{\gminia}(t\del_t)(E_{\gminia})
\subset E_{\gminia}$
for $\nabla^{\reg}_{\gminia}:= \nabla_{\gminia}-d\gminia$. 

Let $F_{0,\gminia}$ denote the endomorphism
of $E_{\gminia|t=0}$
induced by $\nabla^{\reg}_{\gminia}(t\del_t)$.
Because of 
$\bigl[\nabla^{\reg}_{\gminia}(\del_s),
 \nabla^{\reg}_{\gminia}(t\del_t)\bigr]=0$,
the eigenvalues of $F_{0,\gminia}$ are constant.
Since $E_{\gbigk}$ is a strict Deligne-Malgrange lattice,
we have 
$\alpha-\beta\not\in\seisuu-\{0\}$
for any two distinct eigenvalues $\alpha,\beta$
of $F_{0,\gminia}$.

Let us show 
$\nabla^{\reg}_{\gminia}(\del_s)E'_{\gminia}
 \subset E'_{\gminia}$.
Let $\vecv'_{\gminia}$ be a frame of  $E'_{\gminia}$.
Let $A$ and $B$ be determined by the following:
\[
 \nabla^{\reg}_{\gminia}(t\del_t)\vecv'_{\gminia}
=\vecv'_{\gminia}\cdot A,
\quad
 A=\sum_{j=0}^{\infty}A_j\cdot t^j,
\quad
 \nabla^{\reg}_{\gminia}(\del_s)\vecv'_{\gminia}
=\vecv'_{\gminia}\cdot B,
\quad
 B=\sum_{j=-N}^{\infty}B_j\cdot t^j
\]
From the commutativity
$\bigl[\nabla^{\reg}_{\gminia}(t\del_t),\,
 \nabla^{\reg}_{\gminia}(\del_s)\bigr]=0$,
we have the following equation:
\[
 AB+t\del_tB-BA+\del_sA=0.
\]
Assume $N>0$.
Then, we have the equation
$A_0B_{-N}-B_{-N}A_0-NB_{-N}=0$.
Because of $\alpha-\beta\not\in\seisuu-\{0\}$
for two distinct eigenvalues of $A_0$,
we obtain $B_{-N}=0$.
Hence, we have $N\leq 0$,
i.e.,
the entries of $B$ are contained in
$\cnum[\![t,s]\!]$.
\hfill\qed

\subsection{Adjustment of the residue of a logarithmic connection}

Let $k$ be a field whose characteristic number is $0$.
Let $E$ be a free module over $k[\![t]\!]$
with a meromorphic connection
$\nabla$ such that 
$t\nabla(\del_t)(E)\subset E$.
Let $E_0$ denote the specialization of $E$
at $t=0$.
We have the well defined endomorphism
$\Res(\nabla)$ of $E_{0}$.
To distinguish the dependence on $E$,
we denote it by $\Res_E(\nabla)$.
We recall the following standard lemma.

\begin{lem}
\label{lem;07.12.25.1}
We can take a lattice $E'$
of $E\otimes k(\!(t)\!)$
such that (i) $\nabla$ is logarithmic 
with respect to $E'$,
(ii) $\alpha-\beta\not\in\seisuu$
for any distinct eigenvalues of
$\Res_{E'}(\nabla)$.
\end{lem}
\pf
We give only an outline.
Let $S$ denote the set
of the eigenvalues of $\Res_E(\nabla)$.
We say $\alpha<\beta$ for 
$\alpha,\beta\in S$
if $\beta-\alpha\in\seisuu_{>0}$.
We say $\alpha\leq \beta$
if $\alpha=\beta$ or $\alpha<\beta$.
It determines the partial order on $S$.
We put
$\rho(E):=\max\bigl\{
 \beta-\alpha\,\big|\,
 \alpha\leq\beta,\,\,\alpha,\beta\in S
 \bigr\}$.
If $\rho(E)=0$, we have nothing to do.
We will reduce the number $\rho(E)$
by replacing $E$.

Let $S_0$ denote the maximal elements $\beta$ of
$S$ such that there exists $\alpha\in S$ 
with $\alpha<\beta$.
Let $\overline{k}$ denote the algebraically closure
of $k$.
We have the generalized eigen decomposition
$E_{0}\otimes \overline{k}
=\bigoplus_{\alpha\in S} \EE_{\alpha}$.
Note that $S_0$ is preserved
by the action of the Galois group of $\overline{k}$
over $k$.
It is easy to see that 
$\bigoplus_{\alpha\in S_0}\EE_{\alpha}$
comes from the subspace $V$ of $E_{0}$.
Let $E^{(1)}:=t^{-1}\cdot E$.
The specialization $E^{(1)}_{0}$
of $E^{(1)}$ at $t=0$
is naturally isomorphic to $E_{0}$
up to constant multiplication.
Hence, $V$ determines the subspace
$V^{(1)}\subset E^{(1)}_0$.
Let $E^{(2)}$ denote the kernel
of the naturally defined morphism
$E^{(1)}\lrarr E^{(1)}_0/V^{(1)}$.
Then, it can be checked
$\rho(E^{(2)})\leq \rho(E)-1$.
\hfill\qed

\section{Mildly ramified connection}

\subsection{Positive characteristic case}

Let $k$ be an algebraically closed field
whose characteristic number $p$ is positive.
Let $C$ be a smooth divisor of $\Spec^fk[\![s,t]\!]$,
which intersects with the divisor $\{t=0\}$
transversally.
We can take a morphism
$\Spec^fk[\![u]\!]\simeq C\subset \Spec^f k[\![s,t]\!]$
given by $\bigl(s(u),t(u)\bigr)$.
We may assume $t(u)=u$
and $s(u)=u\cdot h(u)$.
We put $s'=s-h(t)\cdot t$.
Then, $C$ is given by the ideal generated by $s'$.
We also have $k[\![s',t]\!]\simeq k[\![s,t]\!]$.
For any positive integer $d$,
we use the notation $s_d'$ to denote a $d$-th root of $s'$.

Let $\nbige$ be a free $k[\![s]\!](\!(t)\!)$-module.
Let $\nabla:\nbige\lrarr \nbige
 \otimes\Omega^1_{k[\![s]\!](\!(t)\!)/k}$
be a flat meromorphic connection.
Let $\psi$ denote the $p$-curvature of $\nabla$.
Let $\Fr$ denote the absolute Frobenius map
$k[\![s]\!](\!(t)\!)\lrarr k[\![s]\!](\!(t)\!)$.
It induces the ring homomorphism
$\Fr:k[\![s]\!](\!(t)\!)[T]\lrarr k[\![s]\!](\!(t)\!)[T]$
by $\Fr\bigl(\sum a_j\cdot T^j\bigr)
  =\sum \Fr(a_j)\cdot T^j$.
Due to the observation of Bost-Laszlo-Pauly
(\cite{laszlo-pauly}.
See also Lemma \ref{lem;07.4.14.2} below),
we have
$P_s(T),P_t(T)\in k[\![s]\!](\!(t)\!)[T]$ such that
$ \det\bigl(T-\psi(\del_s)\bigr)
=\Fr \bigl(P_s\bigr)(T)$ and
$\det\bigl(T-\psi(t\del_t)\bigr)
=\Fr\bigl( P_t\bigr)(T)$.
In general,
the polynomials $P_s(T)$ and $P_t(T)$
have the roots in
$k(\!(s_d)\!)(\!(t_d)\!)$ for some appropriate integer $d$.

\begin{df} \mbox{}
\label{df;07.4.13.1}
We say that $(\nbige,\nabla)$ is mildly ramified at
$\{t=0\}\cup C$,
if the following conditions are satisfied:
\begin{enumerate}
\item
 The roots of the polynomials $P_s(T)=0$ and $P_t(T)=0$
 are contained in  $k[\![s'_d]\!](\!(t_d)\!)$
 for some $d\in\seisuu_{>0}$,
 where $s_d'$ is taken for $C$ as above.
\item
 The roots are of the form $\alpha+\beta$,
 where
 $\alpha\in k[\![s']\!](\!(t_d)\!)$ and
 $\beta\in k[\![s'_d,t_d]\!]$.
\end{enumerate}
We say that $(\nbige,\nabla)$ is mildly ramified,
if it is mildly ramified at $\{t=0\}\cup C$
for some $C$.
\hfill\qed
\end{df}

The connection $\nabla$ induces the relative
connection
$\nabla_t:\nbige\lrarr \nbige\otimes
 \Omega^1_{k[\![s]\!](\!(t)\!)/k[\![s]\!]}$.
We put $\gbigk:=k(\!(s)\!)(\!(t)\!)$
and $\gminik:=k(\!(t)\!)$.
Both of them are equipped with the differential $\del_t$.
We have the natural inclusion
$k[\![s]\!](\!(t)\!)\subset \gbigk$,
and the specialization 
$k[\![s]\!](\!(t)\!)\lrarr \gminik$ at $s=0$.
The morphisms are equivariant with respect to $\del_t$.
Therefore, we have the induced connections
of $\nbige\otimes\gbigk$
and $\nbige\otimes\gminik$,
which are also denoted by $\nabla_t$.

\begin{lem}
\label{lem;07.4.14.10}
Assume that $(\nbige,\nabla)$ is mildly ramified
at $\{t=0\}\cup C$.
Then, the irregular values of
$\bigl(\nbige\otimes\gbigk,\nabla_{t}\bigr)$
are contained in $k[\![s]\!](\!(t_d)\!)_{<0}$,
and their specialization at $s=0$ give
the irregular values for
$\bigl(\nbige\otimes\gminik,\nabla_{t}\bigr)$.
The induced map 
$\Irr(\nbige\otimes\gbigk,\nabla_t)
\lrarr \Irr(\nbige\otimes\gminik,\nabla_t)$
is surjective.
\end{lem}
\pf
Let $\Sol(P_t)$ denote the set of the solutions of $P_t(T)=0$.
By assumption,
any element of $\Sol(P_t)$ is
of the form $\alpha+\beta$ as above.
We have the natural map
$\kappa_1: k[\![s'_d]\!](\!(t_d)\!)\lrarr
 k(\!(s'_d)\!)(\!(t_d)\!)\simeq
 k(\!(s_d)\!)(\!(t_d)\!)$.
The image of $\Sol(P_t)$ via $\kappa_1$
gives the set of the solutions of $P_t(T)=0$
in $k(\!(s_d)\!)(\!(t_d)\!)$.
We remark that the image of $k[\![s'_d,t_d]\!]$ via $\kappa_1$
is contained in $k(\!(s_d)\!)[\![t_d]\!]$.
Hence, we have 
$\kappa_1(\alpha+\beta)_-=\kappa_1(\alpha_-)
 \in k[\![s]\!](\!(t_d)\!)_{<0}$
for any $\alpha+\beta\in \Sol(P_t)$.
Then, the first claim follows from
the characterization of the irregular value
given in Lemma \ref{lem;07.4.14.3}.

On the other hand,
let us take the specialization of $P_{t}(T)$ to $s=0$,
which are denoted by $P_{t,0}(T)\in\gminik[T]$.
Let $\Sol(P_{t,0})$ denote the solution of 
the equation $P_{t,0}(T)=0$,
which is contained in $k(\!(t_d)\!)$ for some appropriate $d$.
Then, $\Sol(P_{t,0})$ is the image of
$\Sol(P_t)$
by the composite $\kappa_2$ of the following morphisms:
\[
 k[\![s_d']\!](\!(t_d)\!)
\simeq
 k[\![s]\!](\!(t_d)\!)[U]\big/\bigl(U^d-s'(s,t)\bigr)
\lrarr
 k(\!(t_d)\!)[U]\big/\bigl(U^d-s'(0,t)\bigr)
\lrarr
 k(\!(t_d)\!)
\]
The last map is given by 
the substitution $U=s'(0,t)^{1/d}\in k(\!(t_d)\!)$
for some choice of $s'(0,t)^{1/d}$.
Any element of 
$k[\![s_d',t_d]\!]$ is mapped into $k[\![t_d]\!]$
via $\kappa_2$,
and the image of any element of
$k[\![s]\!](\!(t_d)\!)=k[\![s']\!](\!(t_d)\!)$ via $\kappa_2$
is given by the natural specialization at $s=0$.
Hence, for any $\kappa_2(\alpha+\beta)\in
 \Sol\bigl(P_{t,0}(T)\bigr)$,
we have 
$\kappa_2(\alpha+\beta)_-=\kappa_2(\alpha_-)$.
Then, the second and third claims follow from 
the characterization of the irregular values
in Lemma \ref{lem;07.4.14.3}.
\hfill\qed

\vspace{.1in}

Let $\varphi:\Spec^f k[\![v]\!]\lrarr \Spec^fk[\![s,t]\!]$
be a morphism given by 
$\varphi^{\ast}(s)=v\cdot\huebar_0(v)$
and $\varphi^{\ast}(t)=v^a$ for some $a>0$.
We assume $a$ is sufficiently smaller than $p$.
We consider the morphism
$\Phi:\Spec^fk[\![u,v]\!]\lrarr \Spec^fk[\![s,t]\!]$
given by
$\Phi^{\ast}s=u+v\cdot \huebar_0(v)$
and $\Phi^{\ast}t=v^a$.
Then, we have
$\Phi^{\ast}s'=
\Phi^{\ast}\bigl(s-h(t)\cdot t\bigr)
=u+v\cdot\huebar_0(v)-h(v^a)\cdot v^a
=u+v\cdot h_1(v)$.
In particular, the divisor
$\Phi^{\ast}(s')=0$ is smooth and
transversal to the divisor $\{v=0\}$.

\begin{lem}
\label{lem;07.4.14.20}
$\Phi^{\ast}(E,\nabla)$ is mildly ramified
at $\{v=0\}\cup \{\Phi^{\ast}(s')=0\}$.
\end{lem}
\pf
Let $\Phi^{\ast}(ds)=a_{1,1}\cdot du+a_{1,2}\cdot dv/v$
and $\Phi^{\ast}(dt/t)=a_{2,1}\cdot du+a_{2,2}\cdot dv/v$,
where $a_{i,j}$ are contained in $k[\![u,v]\!]$.
Then, due to a formula of O. Gabber
(Appendix of \cite{katz3})
we have the following:
\[
 \Phi^{\ast}(\psi)(\del_u)
=a_{1,1}^p\cdot\Phi^{\ast}\bigl(\psi(\del_s)\bigr)
+a_{2,1}^p\cdot\Phi^{\ast}\bigl(\psi(t\del_t)\bigr),
\quad
 \Phi^{\ast}(\psi)(v\del_v)
=a_{1,2}^p\cdot\Phi^{\ast}\bigl(\psi(\del_s)\bigr)
+a_{2,2}^p\cdot\Phi^{\ast}\bigl(\psi(t\del_t)\bigr)
\]
Then, it is easy to check the claim of the lemma
because of the commutativity of
$\psi(\del_s)$ and $\psi(t\del_t)$.
\hfill\qed

\subsection{Mixed characteristic case}

Let $R$ be a subring of $\cnum$
finitely generated over $\seisuu$.
Let $\nbige_R$ be a free $R[\![s]\!](\!(t)\!)$-module,
and let $\nabla:\nbige_R\lrarr
 \nbige_R\otimes\Omega^1_{R[\![s]\!](\!(t)\!)/R}$
be a meromorphic flat connection.
For each $\eta\in S(R)$,
we put $\nbige_{\etabar}:=\nbige_R\otimes_Rk(\etabar)$,
and we have the induced meromorphic flat connection
$\nabla$ of $\nbige_{\etabar}$.

\begin{df}
\label{df;07.4.13.3}
We say that $(\nbige_R,\nabla)$ is 
mildly ramified,
if $(\nbige_{\etabar},\nabla)$ is
mildly ramified for any $\eta\in S(R)$.
Note that the ramification curves may depend 
on $\eta$.
\hfill\qed
\end{df}

If $(\nbige_R,\nabla)$ is mildly ramified,
it is easy to show that
$(\nbige_R,\nabla)\otimes R'$ is also mildly ramified
for any $R'\subset \cnum$ finitely generated over $R$.

\subsection{Complex number field case}

Let $\nbige_{\cnum}$ be a free
$\cnum[\![s]\!](\!(t)\!)$-module
with a meromorphic connection
$\nabla:\nbige_{\cnum}\lrarr\nbige_{\cnum}\otimes
 \Omega^1_{\cnum[\![s]\!](\!(t)\!)/\cnum}$.

\begin{df}
We say that $(\nbige_{\cnum},\nabla)$ is algebraic,
if there exists a subring $R\subset\cnum$
finitely generated over $\seisuu$,
a free $R[\![s]\!](\!(t)\!)$-module $\nbige_R$
with a meromorphic connection
$\nabla:\nbige_R\lrarr
 \nbige_R\otimes\Omega^1_{R[\![s]\!](\!(t)\!)/R}$
such that
$(\nbige_R,\nabla)\otimes_R\cnum
\simeq
 (\nbige_{\cnum},\nabla)$.
Such $(\nbige_R,\nabla)$ is called
an $R$-model of $(\nbige_{\cnum},\nabla)$.
\hfill\qed
\end{df}

\begin{df}
Let $(\nbige_{\cnum},\nabla)$ be algebraic.
We say $(\nbige_{\cnum},\nabla)$ is mildly ramified,
if an $R$-model of $(\nbige_{\cnum},\nabla)$
is mildly ramified for some $R$.
\hfill\qed
\end{df}

We put $\gbigk_{\cnum}:=\cnum(\!(s)\!)(\!(t)\!)$
and $\gminik_{\cnum}:=\cnum(\!(t)\!)$.
We have the induced relative connection
$\nabla_t:\nbige_{\cnum}\lrarr\nbige_{\cnum}
 \otimes\Omega^1_{\cnum[\![s]\!](\!(t)\!)/\cnum[\![s]\!]}$.
We put
$\bigl(\nbige_{\gbigk_{\cnum}},\nabla_{t}\bigr)
 :=\bigl(\nbige_{\cnum},\nabla_t\bigr)\otimes\gbigk_{\cnum}$
and $\bigl(\nbige_{\gminik_{\cnum}},\nabla_t\bigr):=
 \bigl(\nbige_{\cnum},\nabla_t\bigr)\otimes\gminik_{\cnum}$.

\begin{prop}
\label{prop;07.4.14.15}
Assume that $(\nbige_{\cnum},\nabla)$
is algebraic and mildly ramified.
Then the irregular values of
$(\nbige_{\gbigk_{\cnum}},\nabla_t)$ are 
contained in $\cnum[\![s]\!](\!(t_d)\!)_{<0}$
for some $d\in\seisuu_{>0}$,
and their specializations at $s=0$
give the irregular values of
$(\nbige_{\gminik_{\cnum}},\nabla_t)$.
The induced map
$\Irr(\nbige_{\gbigk_{\cnum}},\nabla_t)
\lrarr
 \Irr(\nbige_{\gminik_{\cnum}},\nabla_t)$
is surjective.
\end{prop}
\pf
We take a subring $R\subset\cnum$ 
finitely generated over $\seisuu$,
and an $R$-model $(\nbige_R,\nabla)$ such that
$(\nbige_R,\nabla)\otimes \cnum
 \simeq (\nbige_{\cnum},\nabla)$.
We may assume that the irregular decomposition
of $(\nbige_{\gbigk_{\cnum}},\nabla_t)$
is defined on $R(\!(s_d)\!)(\!(t_d)\!)$
(Corollary \ref{cor;07.4.28.1}):
\begin{equation}
 \label{eq;07.4.13.2}
 (\nbige_{R},\nabla_t)
 \otimes R(\!(s_d)\!)(\!(t_d)\!)
=\bigoplus_{\gminia\in\Irr(\nbige_{\gbigk_{\cnum}},\nabla_t)}
 \bigl(\nbige_{\gminia},\nabla_{\gminia,t}\bigr)
\end{equation}
We may also have a Deligne-Malgrange lattice
$\bigoplus E_{\gminia}\subset\bigoplus\nbige_{\gminia}$.

Let $p$ be a sufficiently large prime,
and let $\eta$ be any point of $S(R,p)$.
We put $\gbigk_{\etabar}:=k(\etabar)(\!(s_d)\!)(\!(t_d)\!)$
and $\gminik_{\etabar}:=k(\etabar)(\!(t_d)\!)$.
We have the decomposition
of $(\nbige_{\gbigk_{\etabar}},\nabla_t):=
 (\nbige_R,\nabla_t)\otimes \gbigk_{\etabar}$
induced by (\ref{eq;07.4.13.2}):
\[
 (\nbige_{\gbigk_{\etabar}},\nabla_t)
=\bigoplus_{\gminia\in\Irr(\nbige_{\gbigk_{\cnum}},\nabla_t)}
 \bigl(\nbige_{\gminia,\etabar},\nabla_{\gminia,t}\bigr)
\]
We use the notation $\nbigf_{\etabar}$
to denote the naturally induced morphism
$R(\!(s_d)\!)(\!(t_d)\!)
\lrarr \gbigk_{\etabar}$
and $R(\!(t_d)\!)\lrarr \gminik_{\etabar}$.
Since $\nabla_{\gminia,t}
-d\gminia\cdot \id_{\nbige_{\gminia,\etabar}}$
are logarithmic with respect to the lattice
$E_{\gminia,\etabar}$,
we can conclude that
$\Irr\bigl(\nbige_{\gbigk_{\etabar}},\nabla_t\bigr)$
is the image of 
$\Irr\bigl(\nbige_{\gbigk_{\cnum}},\nabla_t\bigr)$
via the map $\nbigf_{\etabar}$.
Due to Lemma \ref{lem;07.4.14.10},
$\nbigf_{\etabar}(\gminia)$ are contained in
$k(\etabar)[\![s]\!](\!(t_d)\!)_{<0}$
for any $\gminia\in \Irr(\nbige_{\gminik_{\cnum}},\nabla_t)$.
Then, it follows that
$\gminia$ are contained in $R[\![s]\!](\!(t_d)\!)_{<0}$.

Moreover,
$\nbigf_{\etabar}\bigl(\gminia_{|s=0}\bigr)
=\nbigf_{\etabar}\bigl(\gminia\bigr)_{|s=0}$
give the irregular values of
$(\nbige_{\gminik_{\etabar}},\nabla_t)$
for any $\gminia\in\Irr(\nbige_{\gbigk_{\cnum}},\nabla_t)$,
due to Lemma \ref{lem;07.4.14.10}.
To conclude that
$\gminia_{|s=0}$ gives the irregular values of
$(\nbige_{\gminik_{\cnum}},\nabla_t)$,
we use the following lemma.

\begin{lem}
\label{lem;07.4.14.11}
Let $(\nbige,\nabla)$ be a meromorphic connection on $R(\!(t)\!)$.
Let $\gminia\in R(\!(t_d)\!)_{<0}$.
If $\nbigf_{\etabar}(\gminia)$ are the irregular values
for $(\nbige,\nabla)_{\etabar}$ on $k(\etabar)(\!(t)\!)$
for any $\etabar$,
then $\gminia$ is an irregular value
for $(\nbige,\nabla)\otimes\cnum(\!(t)\!)$.
\end{lem}
\pf
We may assume to have
the irregular decomposition
$(\nbige,\nabla)=\bigoplus_i
 (\nbige_{i},d\gminia_i+\nabla^{\reg}_{i})$
on $R(\!(t_d)\!)$,
due to Corollary \ref{cor;07.4.28.1}.
Then, for some $i$,
there are infinitely many $\eta\in S(R)$ such that
$\nbigf_{\etabar}(\gminia)-\nbigf_{\etabar}(\gminia_i)=0$
in $k(\etabar)(\!(t_d)\!)_{<0}$.
It implies $\gminia=\gminia_i$.
Thus, we obtain Lemma \ref{lem;07.4.14.11}.
\hfill\qed

\vspace{.1in}

Let us return to the proof of
Proposition \ref{prop;07.4.14.15}.
Let $\gminib\in\Irr(\nbige_{\gminik_{\cnum}},\nabla_t)$.
Because of the surjectivity in Lemma \ref{lem;07.4.14.10},
there exists
$\gminia\in \Irr(\nbige_{\gbigk_{\cnum}},\nabla_t)$
such that 
$\nbigf_{\etabar}(\gminia_{|s=0})
=\nbigf_{\etabar}(\gminib)$ in 
$k(\!(\etabar)\!)(\!(t_d)\!)_{<0}$
for infinitely many $\eta\in S(R)$.
It implies $\gminia_{|s=0}=\gminib$.
Hence, we obtain the surjectivity of 
the induced map
$\Irr(\nbige_{\gbigk_{\cnum}},\nabla_t)
\lrarr
 \Irr(\nbige_{\gminik_{\cnum}},\nabla_t)$.
Thus the proof of Proposition \ref{prop;07.4.14.15}
is finished.
\hfill\qed

\vspace{.1in}

Let $\varphi_{\cnum}:\Spec^f\cnum[\![v]\!]
 \lrarr \Spec^f\cnum[\![s,t]\!]$
be an algebraic morphism,
i.e.,
there exist a morphism $\Spec A_1\lrarr \Spec A_2$
for some regular rings $A_i$ $(i=1,2)$ finitely generated over $\cnum$,
such that the completion at some closed points
is isomorphic to $\varphi_{\cnum}$.
We assume $\varphi_{\cnum}^{\ast}(t)\neq 0$.
We have the induced map
$\varphi^{\ast}_{<0}:
 \cnum[\![s]\!](\!(t_d)\!)\big/\cnum[\![s,t_d]\!]
\lrarr \cnum(\!(v_d)\!)\big/\cnum[\![v_d]\!]$
for any $d$.

\begin{prop}
\label{prop;07.6.19.1}
If $(\nbige_{\cnum},\nabla)$ is algebraic and mildly ramified,
the set of the irregular values of
$\varphi_{\cnum}^{\ast}(\nbige_{\cnum},\nabla)$
is given by the image of
$\Irr\bigl(\nbige_{\cnum},\nabla_t\bigr)$
via $\varphi^{\ast}_{<0}$.
\end{prop}
\pf
By extending $R$,
we may assume that $\varphi$ is induced from
$\varphi_R:\Spec^f R[\![v]\!]\lrarr \Spec^f R[\![s,t]\!]$
given by 
$\varphi^{\ast}_R(t)=v^a$ and 
$\varphi_R^{\ast}(s)=v\cdot h(v)$.
We have the induced map
$\varphi_R^{\ast}: R(\!(s_d)\!)(\!(t_d)\!)
 \lrarr R(\!(v_d)\!)$.
Let $\Phi:\Spec^f R[\![u,v]\!]\lrarr \Spec^f R[\![s,t]\!]$
be given by $t=v^a$ and $s=u+v\cdot h(v)$.
Then, $\Phi^{\ast}(\nbige,\nabla)$ is mildly ramified
because of Lemma \ref{lem;07.4.14.20}.

We put $\gbigk(u,v):=\cnum(\!(u_d)\!)(\!(v_d)\!)$
and
$\nbige_{\gbigk(u,v)}:=\nbige_{\cnum}\otimes\gbigk(u,v)$
on which the relative connection $\nabla_v$ is induced.
We have the induced map
$R[\![s]\!](\!(t_d)\!)\big/R[\![s,t_d]\!]
 \lrarr
 R[\![u]\!](\!(v_d)\!)\big/
 R[\![u,v_d]\!]$,
which is denoted by $\Phi^{\ast}_{<0}$.
Then, we have only to show that
$\Irr(\nbige_{\gbigk(u,v)},\nabla_v)$
is the same as the image of
$\Irr\bigl(\nbige_{\gbigk_{\cnum}},\nabla_t\bigr)$
via the map $\Phi^{\ast}_{<0}$
due to Proposition \ref{prop;07.4.14.15}.
Since both of them are contained in
$\cnum[\![u]\!](\!(v_d)\!)/\cnum[\![u,v_d]\!]$,
we have only to compare them in 
$\cnum(\!(u)\!)(\!(v_d)\!)/\cnum(\!(u)\!)[\![v_d]\!]$.

The meromorphic connection
$\bigl(
 \nbige_{\gbigk_{\cnum}}\otimes
 \cnum(\!(s)\!)(\!(t_d)\!),
 \nabla_t
\bigr)$ is unramified,
because the irregular values are contained
in $\cnum[\![s]\!](\!(t_d)\!)
 \big/\cnum[\![s,t_d]\!]$.
By Lemma \ref{lem;07.12.25.1},
we have a strict Deligne-Malgrange lattice $E_{\cnum}$
which is the free $\cnum(\!(s)\!)[\![t_d]\!]$-module,
and 
the irregular decomposition
with respect to the relative connection $\nabla_t$:
\[
 \bigl(
 E_{\cnum},
\nabla_t
\bigr)
=\bigoplus_{\gminia\in\Irr(\nbige_{\gbigk_{\cnum}},\nabla_t)}
 \bigl(E_{\gminia},\nabla_{\gminia,t}\bigr)
\]
Due to the uniqueness of the irregular decomposition
and the commutativity of
$\nabla(\del_s)$ and $\nabla(t\del_t)$,
it is standard to show that
$\nabla(\del_s)\bigl(E_{\gminia}(\!(t_d)\!)\bigr)
\subset E_{\gminia}(\!(t_d)\!)$.
(See the proof of Lemma \ref{lem;07.9.19.14},
for example.)
Hence, it is the decomposition of the meromorphic flat connection:
\[
\bigl(E_{\cnum},\nabla\bigr)
=\bigoplus \bigl(E_{\gminia},\nabla_{\gminia}\bigr).
\]
We put $\nabla'_{\gminia}:=\nabla_{\gminia}-d\gminia$.
By construction,
we have $\nabla'_{\gminia}(t_d\del_{t_d})(E_{\gminia})
\subset E_{\gminia}$.
Since $E_{\cnum}$ is assumed to be strict Deligne-Malgrange,
it can be shown that
$\nabla'_{\gminia}(\del_s) (E_{\gminia})
 \subset E_{\gminia}$
by a standard argument.
(See the last part of the proof of Proposition \ref{prop;07.9.20.1},
for example.)
We put $\nabla'=\bigoplus \nabla'_{\gminia}$.

Let $\vecv$ be a frame of $E_{\cnum}$
compatible with the irregular decomposition.
Let $A$ and $B$ be determined by
$\nabla'\vecv=\vecv\cdot (A\cdot dt_{d}/t_d+B\cdot ds)$.
Then, $A,B\in M_r\bigl(\cnum(\!(s)\!)[\![t_d]\!]\bigr)$.
We remark 
$\Phi^{\ast}(s)^{-k}\in \cnum(\!(u)\!)[\![v]\!]$
for any integer $k$.
Then, it is easy to see that
$E_{\cnum}\otimes \cnum(\!(u)\!)[\![v_d]\!]$
gives a Deligne-Malgrange lattice
of $\nbige\otimes\cnum(\!(u)\!)(\!(v_d)\!)$
with respect to $\nabla(v_d\del_{v_d})$,
and the irregular decomposition of $\Phi^{\ast}(\nbige,\nabla)$
is given as follows:
\[
 \nbige\otimes\cnum(\!(u)\!)(\!(v_d)\!)
\simeq
 \bigoplus_{\gminib\in
 \cnum(\!(u)\!)(\!(v_d)\!)/\cnum(\!(u)\!)[\![v_d]\!]}
 \left(
 \bigoplus_{\Phi^{\ast}_{<0}(\gminia)=\gminib}
 E_{\gminia}\otimes\cnum(\!(u)\!)(\!(v_d)\!)
 \right).
\]
Thus, we are done.
\hfill\qed

\section{Resolution of turning points}
\label{section;07.4.28.2}

\subsection{Resolution of the discriminants of polynomials}

Let $R$ be a regular subring of $\cnum$
which is finitely generated over $\seisuu$.
Let $X_{R}$ be a smooth projective surface over $R$.
Let $D_R$ be a simply effective normal crossing divisor
of $X_R$.
We assume that $X_R\otimes_{\seisuu}\seisuu/p\seisuu$
is smooth or empty for each $p$.
Let $N$ be a positive integer.

Take $\eta\in S(R,p)$.
We put $X_{\etabar}:=X_R\otimes_Rk(\etabar)$.
We denote the function field of $X_{\etabar}$
by $K(X_{\etabar})$.
Let $\nbigp^{(a)}(T)\in
 \bigoplus_{j=0}^r H^0\bigl(X_{\etabar},
 \nbigo_X(jND_{\etabar})\bigr)\cdot T^{r-j}$
be monic polynomials
$(a=1,\ldots,L)$.
The tuple $\bigl(\nbigp^{(a)}\,\big|\,a=1,\ldots,L\bigr)$
is denoted by $\boldnbigp$.
We regard them as elements of $K(X_{\etabar})[T]$.
Let $\nbigp^{(a)}=\prod_{i=1}^{m(a)} (\nbigp_i^{(a)})^{e(i,a)}$
be the irreducible decomposition.
The monic polynomials
$\nbigp_i^{(a)}$ are contained in
$\bigoplus_{j=0}^{r_i(a)}
 H^0\bigl(X_{\etabar},\nbigo_X(jND_{\etabar})\bigr)
 \cdot T^{r_i(a)-j}$,
where $r_i(a):=\deg_T\nbigp_i^{(a)}$.
We regard the discriminants
$\disc(\nbigp_i^{(a)})$ as the elements of 
the function field $K(X_{\etabar})$.
There exists a constant $M_1>0$,
which is independent of the choice of $\etabar$ and $p$,
such that 
$\disc(\nbigp_i^{(a)})$ are contained in
$H^0\bigl(X_{\etabar},\nbigo_X(M_1\cdot D_{\etabar})\bigr)$.
We put as follows:
\[
 \disc(\boldnbigp):=
 \prod_{a=1}^L
 \prod_{j=1}^{m(a)} \disc(\nbigp_j^{(a)})
\in K(X_{\etabar})
\]
There exists a constant $M_2>0$,
which is independent of the choice of $\etabar$ and $p$,
such that $\disc(\boldnbigp)$ is contained in
$H^0\bigl(X_{\etabar},\nbigo_X(M_2\cdot D_{\etabar})\bigr)$.
Let $Z(\boldnbigp)$ denote the $0$-set of
$\disc(\boldnbigp)$,
when we regard $\disc(\boldnbigp)$ is a section of
the line bundle $\nbigo_X(M_2\cdot D_{\etabar})$.
We may assume $D_{\etabar}\subset Z(\boldnbigp)$,
by making $M_2$ larger.
Since $Z(\boldnbigp)$ is a member
of some bounded family,
the following lemma immediately follows from
the flattening lemma (see \cite{mumford})
and the semi-continuity theorem (see \cite{har1})
for the flat family.

\begin{lem}
There exists a constant $M_3$,
which is independent of $\etabar$ and $p$,
such that
the arithmetic genus of $Z(\boldnbigp)$
is smaller than $M_3$.
\hfill\qed
\end{lem}

Let $P$ be any closed point of $X_{\etabar}$.
We put 
$(X^{(0)}_{\etabar},P^{(0)})
 :=(X_{\etabar},P)$.
Inductively,
let $\pi^{(i)}:
 X^{(i)}_{\etabar}\lrarr X^{(i-1)}_{\etabar}$
be the blow up at $P^{(i-1)}$,
and let us take a point 
$P^{(i)}\in \pi^{(i)}(P^{(i-1)})$.
Let $\pi_i$ denote the naturally induced map
$X^{(i)}\lrarr X$.
By the classical arguments
(see Section V.3 in \cite{har1}, for example),
we can show the following lemma.
\begin{lem}
\label{lem;07.4.15.1}
There exists some $i_0$, independent of the choice of
$p$, $\etabar$ and the points $P^{(i)}$,
such that the divisor 
 $(\pi^{(i)})^{-1}Z(\boldnbigp)$
is normal crossing around 
the exceptional divisor $(\pi^{(i)})^{-1}(P^{(i-1)})$
for any $i\geq i_0$.
\end{lem}
\pf
We give only an outline.
We use the notation $p_a$ to denote the arithmetic genus.
Let $Y$ denote the reduced scheme
associated to $Z(\boldnbigp)$.
Let $Y=\bigcup Y_j$ denote the irreducible decomposition.
We have $p_a(Y)\leq M_3$ and $p_a(Y_j)\leq M_3$.
Let $\Ytilde_i$ denote the inverse image of
$Y$ via $\pi_i$ with the reduced structure.
Let $\Ytilde_{i,j}$ denote the strict transform
of $Y_j$ via $\pi_i$.
Let $C_{i,q}$ denote the strict transform
of $(\pi^{(q)})^{-1}(P^{(q-1)})$
via the natural map
$X^{(i)}_{\etabar}\lrarr X^{(q)}_{\etabar}$.
We have
$\Ytilde_{i}=
 \bigcup_{j}\Ytilde_{i,j}\cup
 \bigcup_{q}C_{i,q}$.
Let $r_{P^{(q)}}(\Ytilde_q)$ denote
the multiplicity of $P^{(q)}$ in $\Ytilde_q$.
We use the notation $r_{P^{(q)}}(\Ytilde_{q,j})$
in a similar meaning.
We have the equality (Section V.3 of \cite{har1}):
\[
  p_a(\Ytilde_{i,j})
=p_a(Y_j)
-\sum_{q\leq i-1}\frac{1}{2}
 r_{P^{(q)}}(\Ytilde_{q,j}) \cdot
 \bigl(r_{P^{(q)}}(\Ytilde_{q,j})-1 \bigr)
\]
By our choice,
$P^{(i)}$ is a smooth point of 
$\Ytilde_{i,j}$ for any  $i\geq i(1)$
if $P^{(i(1))}$ is a smooth point of
$\Ytilde_{i(1),j}$.
Hence, we obtain
$r_{P^{(i)}}(\Ytilde_{i,j})\leq 1$ 
if $i$ is sufficiently large.
We also have the following equality (Section V.3 of \cite{har1}):
\[
 p_a(\Ytilde_i)=p_a(Y)
-\sum_{q\leq i-1}\frac{1}{2}
 \bigl( r_{P^{(q)}}(\Ytilde_{q})-1\bigr)
\cdot
  \bigl( r_{P^{(q)}}(\Ytilde_{q})-2\bigr)
\]
Assume $r_{P^{(q)}}(\Ytilde_q)=2$.
Then, as explained in the proof of
Theorem 3.9 in Section V of \cite{har1},
there are three possibility:
\begin{itemize}
\item
 $\Ytilde_q$ is normal crossing around $P^{(q)}$.
\item
 Let $\Ytilde_{q+1}'$ denote the strict transform 
 of  $\Ytilde_{q}$ via $\pi^{(q+1)}$.
 Then, it is nonsingular
 in a neighbourhood of $(\pi^{(q+1)})^{-1}(P^{(q)})$,
 and $\Ytilde_{q+1}'$ and $(\pi^{(q+1)})^{-1}(P^{(q)})$
 intersect at one point with multiplicity $2$.
 If $P^{(q+1)}$ and $P^{(q+2)}$ are also 
 singular points of $\Ytilde_{q+1}$
 and $\Ytilde_{q+2}$ respectively,
 we have $r_{P^{(q+2)}}(\Ytilde_{q+2})=3$.
\item
 $\Ytilde_{q+1}'$ and 
 $(\pi^{(q+1)})^{-1}(P^{(q)})$
 intersects at one point, 
 whose multiplicity in $\Ytilde_{q+1}'$ is $2$.
 If $P^{(q+1)}$ is singular point of a
 $\Ytilde_{q+1}$,
 we have 
 $r_{P^{(q+1)}}(\Ytilde_{q+1})=3$.
\end{itemize}
Hence, we obtain that
$\Ytilde_i$ are normal crossing around $P^{(i)}$
for sufficiently large $i$.
\hfill\qed

\vspace{.1in}

Let $i\geq i_0$.
Let $C_i(\boldnbigp)$ denote
 the closure of 
 $Z(\boldnbigp)\cap (X_{\etabar}-D_{\etabar})$ 
in $X_{\etabar}^{(i)}$.
We take a local coordinate neighbourhood
$(U^{(i)},s^{(i)},t^{(i)})$ around $P^{(i)}$
such that
(i) $(t^{(i)})^{-1}(0)$ is
 $U^{(i)}\cap \bigl(\pi^{(i)}\bigr)^{-1}(P^{(i-1)})$,
(ii) if $P^{(i)}$ is contained in $C_i(\boldnbigp)$,
   then $(s^{(i)})^{-1}(0)=U^{(i)}\cap C_i(\boldnbigp)$,
(ii)' if $P^{(i)}$ is not contained in $C_i(\boldnbigp)$,
 then $s^{(i)}$ may be anything.
Because of generalized Abhyankar's lemma
(see Expose XIII Section 5 of \cite{raynaud}),
any solutions of the equations
$\pi_i^{\ast}\nbigp_j^{(a)}(T)=0$ 
$(a=1,\ldots,L,\,\,j=1,\ldots,m(a))$
are contained
in $k(\etabar)[\![s_d^{(i)}]\!](\!(t_d^{(i)})\!)$
for some appropriate $d$,
which is a factor of $r!$.

\begin{lem}
\label{lem;07.4.16.1}
There exists an $i_1$, which is independent of the choice of
$\etabar$, $p$ and the points $P^{(i)}$,
such that the following holds
for any $i\geq i_1$:
\begin{itemize}
\item
Any solutions of the equations
$\pi_i^{\ast}\nbigp_j^{(a)}(T)=0$ 
$(a=1,\ldots,L,\,\,j=1,\ldots,m(a))$
are contained in the following:
\[
 k(\etabar)[\![s_d^{(i)},t_d^{(i)}]\!]+
      k(\etabar)[\![s^{(i)}]\!](\!(t_d^{(i)})\!)
\]
\end{itemize}
\end{lem}
\pf
If $P^{(i_0)}$ is not contained in
$C_{i_0}(\boldnbigp)$,
then $P^{(i)}\not\in C_{i}(\boldnbigp)$ for any
$i\geq i_0$,
and the claim is obvious in this case.
Assume $P^{(i_0)}$ is contained in $C_{i_0}(\boldnbigp)$.
Let $\alpha^{(i_0)}_l\in
k(\etabar)[\![s_d^{(i_0)}]\!](\!(t_d^{(i_0)})\!)$
be any solution of $\pi^{\ast}_{i_0}\nbigp_j^{(a)}(T)=0$
for some $(a,j)$.
Note that there exists a constant $M_4$,
which is independent of the choice of
 $\etabar$, $p$, and
the sequence of the points $P^{(i)}$,
with the following property:
\begin{itemize}
\item
The orders of the poles of the coefficients of
$\pi^{\ast}_{i_0}\nbigp_j^{(a)}(T)$
with respect to $t^{(i_0)}$ are dominated by $M_4$,
\end{itemize}
Hence, there exists a constant $M_5$,
which is independent of the choice of $\etabar$, $p$,
the sequence of the points $P^{(i)}$,
$(a,j)$ and $\alpha^{(i_0)}_l$,
with the following property:
\begin{itemize}
\item
The order of the pole of $\alpha^{(i_0)}_l$
with respect to $t_d^{(i_0)}$ are dominated by $M_5$.
\end{itemize}
If $P^{(i)}$ are contained in $C_{i}(\boldnbigp)$
for $i\geq i_0$,
we may assume
$(\pi^{(i)})^{\ast}(s^{(i-1)})=s^{(i)}\cdot t^{(i)}$
and $(\pi^{(i)})^{\ast}(t^{(i-1)})=t^{(i)}$.
Hence, 
the pull back of $\alpha^{(i_0)}_l$
via $X^{(i)}_{\etabar}\lrarr X^{(i_0)}_{\etabar}$
are contained in 
$k[\![s_d^{(i)},t_d^{(i)}]\!]$
for sufficiently large $i$.
\hfill\qed

\subsection{Proof of Theorem \ref{thm;07.4.15.10}}

If we take a sufficiently large $R$,
then
$\nbige$, $\nabla$, $X$ and $D$ come from
$\nbige_R$, $\nabla_R$, $X_R$ and $D_R$
which are defined over $R$.
We may also assume that we have 
the canonical lattice 
$E_R\subset\nbige_R$ defined over $R$.
(See \cite{malgrange})
By applying a theorem of Sabbah (\cite{sabbah4}),
we may assume that any cross points of $D$ are not turning.
Let $P$ be a turning point contained in a smooth part of $D$.
Let $U$ be a neighbourhood of $P$
with an \'{e}tale morphism $(x,y):U\lrarr A^2$
such that $x^{-1}(0)=D\cap U$.
For simplicity,
$U$ does not contain 
any other turning points than $P$.
We may assume $P$ and $(U,x,y)$ are also defined over $R$.
We have only to take a proper birational map
$\pi:U'\lrarr U$ such that 
$\pi^{-1}(\nbige,\nabla)$ has no turning points.
On $U$, we have the vector field
$t\del_t$ and $\del_s$.
By taking blow up of $X$ outside of $U$,
and by extending $D$,
we may assume that $x\del_x$ and $\del_y$
are sections of $\Theta_X(M_0D)$.
We have a positive number $M_0'$
such that $\nabla(E_R)\subset 
 E_R(M_0'D_R)$.
Hence, we have the constant $M_1$
such that $\nabla(x\del_x)(E_R)$
and $\nabla(\del_y)(E_R)$ are contained in
$E_R(M_1D_R)$.

Let $p$ be a large prime.
For $\eta\in S(R,p)$,
let $\nbige_{\etabar},E_{\etabar},
 \nabla_{\etabar},X_{\etabar},D_{\etabar},
 P_{\etabar}$ and $U_{\etabar}$ denote
the induced objects over $k(\etabar)$.
Let $\psi$ be the $p$-curvature of
$(\nbige_{\etabar},\nabla_{\etabar})$.
We put $\psi_x:=\psi(x\del_x)$
and $\psi_y:=\psi(\del_y)$.
Because of 
$\nabla(\del_x)\bigl(E_{\etabar}\bigr)
 \subset E_{\etabar}\bigl(M_1D_{\etabar}\bigr)$
and $\nabla(y\del_y)\bigl(E_{\etabar}\bigr)
 \subset E_{\etabar}(M_1D_{\etabar})$,
we have
$\psi_x,\psi_y\in\End(E_{\etabar})\otimes\nbigo(pM_1D)$.
Hence,
the characteristic polynomials
$ \det(T-\psi_x)$ and
$\det(T-\psi_y)$ are contained in
$ \bigoplus_{j=0}^n
 H^0\bigl(X_{\etabar},\nbigo(pjM_1D)\bigr)\cdot T^{n-j}$.
Due to the excellent observation of
Bost, Laszlo and Pauly \cite{laszlo-pauly},
we have the following lemma.
\begin{lem}
\label{lem;07.4.14.2}
We have $P_x(T)$ and $P_y(T)$
in $\bigoplus _{j=0}^n
 H^0\bigl(X_{\etabar},\nbigo(jM_1D)\bigr)\cdot T^{n-j}$
such that
$\det(T-\psi_x)=\Fr^{\ast}P_x(T)$ and
$\det(T-\psi_y)=\Fr^{\ast}P_y(T)$,
where $\Fr:X_{\etabar}\lrarr X_{\etabar}$
denotes the absolute Frobenius morphism.
\end{lem}
\pf
We reproduce the argument in \cite{laszlo-pauly}
for the convenience of the reader.
Let $T'$ be a formal variable.
Because of the Cartier descent,
we have only to show $\del_y\det(1-T'\psi_{\kappa})=0$
and $x\del_x\det(1-T'\psi_{\kappa})=0$
for $\kappa=x,y$.
Let $\vecv$ be a local frame of $E_{\etabar}$.
Let $A$, $B$, and $\Psi_{\kappa}$ $(\kappa=x,y)$
be determined by
$\nabla\vecv=\vecv(Ady+Bdx/x)$,
$\psi_{\kappa}\vecv=\vecv\cdot \Psi_{\kappa}$.
Because of $(\del_y)^p=0$ and $(x\del_x)^p=x\del_x$,
we have
$\psi_x=\nabla(x\del_x)^p-\nabla(x\del_x)$
and $\psi_y=\nabla(\del_y)^p$.
Since $\nabla$ is flat,
we have the commutativity
$[\nabla(\del_y),\nabla(x\del_x)]
=[\nabla(\del_y),\nabla(\del_y)]
=[\nabla(x\del_x),\nabla(x\del_x)]=0$.
Hence, we have
$[\nabla(\del_y),\psi_{\kappa}]
=[\nabla(x\del_x),\psi_{\kappa}]=0$
for $\kappa=x,y$.
Therefore,
$\del_y\Psi_{\kappa}+[A,\Psi_{\kappa}]
=x\del_x\Psi_{\kappa}+[B,\Psi_{\kappa}]=0$,
and thus
$\Tr(\Psi_{\kappa}^n\del_y\Psi_{\kappa})
=\Tr(\Psi_{\kappa}^nx\del_x\Psi_{\kappa})=0$.

Recall $\del_y\det(M)=\det(M)\cdot \Tr(M^{-1}\del_yM)$
for an invertible matrix $M$.
For $M=\id-T'\Psi_{\kappa}$,
we have
$M^{-1}=\sum T^{\prime\,n}\Psi_{\kappa}^n$
and 
\[
 \del_y\det(\id-T'\Psi_{\kappa})
=-T'\det(\id-T'\Psi_{\kappa})
\cdot\sum_{n=0}T^{\prime\,n}
 \Tr(\Psi_{\kappa}^n\del_y\Psi_{\kappa})
=0
\]
Similarly,
we have $x\del_x\det(\id-T'\Psi_s)=0$.
\hfill\qed

\vspace{.1in}

Inductively, we construct 
the blow up $\pi^{(i)}:X^{(i)}\lrarr X^{(i-1)}$
as follows.
First, let $\pi^{(1)}:X^{(1)}\lrarr X$ be the blow up at $P$,
and we put $\pi_1:=\pi^{(1)}$.
Let $\pi^{(2)}:X^{(2)}\lrarr X^{(1)}$
denote the blow up at the turning points
of $\pi_1^{\ast}(\nbige,\nabla)$
contained in $\pi_1^{-1}(U)$,
and we put $\pi_2:=\pi^{(1)}\circ\pi^{(2)}$.
When $\pi^{(i)}:X^{(i)}\lrarr X^{(i-1)}$ is given,
let $\pi_i:X^{(i)}\lrarr X$  denote the naturally induced morphism,
and let $\pi^{(i+1)}:X^{(i+1)}\lrarr X^{(i)}$ be 
the blow up at the turning points of
$\pi_i^{\ast}(\nbige,\nabla)$ contained in 
$\pi_i^{-1}(U)$.

We take subrings $R^{(i)}\subset \cnum$
such that
(i)  $R^{(i-1)}\subset R^{(i)}$,
 and $R^{(i)}$ is smooth and finitely generated over $R^{(i-1)}$,
(ii) $X^{(i)}$ and the turning points contained in $\pi_i^{-1}(U)$
 are defined over $R^{(i)}$.
Let $\eta\in S(R,p)$.
We take geometric points $\etabar(i)$ of $S(R^{(i)},p)$
for any $i$
with the morphisms
$\etabar(i)\lrarr \etabar(i-1)\lrarr \eta$
compatible with $\Spec R^{(i)}\lrarr \Spec R^{(i-1)}\lrarr \Spec R$.
For $j\leq i$,
$X^{(j)}$ are defined over $R^{(i)}$,
and we have
$X^{(j)}_{R^{(i)}}\otimes_{R^{(i)}}k(\etabar(i))
\simeq
 X^{(j)}_{\etabar(j)}\otimes_{\etabar(j)}k(\etabar(i))$.
And the objects over them are naturally related by the pull backs.

Let $P_{\kappa}(T)=\prod P_{\kappa,j}(T)^{e(\kappa,j)}$
denote the irreducible decomposition of the polynomials
$P_{\kappa}(T)$ above $(\kappa=x,y)$.
Applying Lemma \ref{lem;07.4.15.1}
and Lemma \ref{lem;07.4.16.1},
we can show that
there exist $i_1$ and $p_1$ such that
the following claims hold for any $i\geq i_1$, $p\geq p_1$
and any $\etabar(i)$:
\begin{itemize}
\item
 Let $C$ be any exceptional divisor
 with respect to $\pi^{(i)}_{\etabar(i)}$.
 Then, 
 $\bigcup_{j,\kappa}\pi_{i,\etabar(i)}^{-1}
 \bigl(\disc(P_{\kappa,j})\cup D\bigr)$
 are normal crossing around $C$.
\item
 Let $Z_{\kappa,j}^{(i)}$ denote the closure of
 $\disc(P_{\kappa,j})\cap (X_{\etabar(i)}-D_{\etabar(i)})$
 in $X^{(i)}_{\etabar(i)}$.
 If $C$ intersects at $Q$ with
 $Z_{\kappa,j}^{(i)}$ for some $(\kappa,j)$,
 we take a coordinate neighbourhood
 $(U_Q,z,w)$ such that
 $w^{-1}(0)=C\cap U_Q$
 and $z^{-1}(0)$ is $Z_{\kappa,j}^{(i)}\cap U_Q$.
 Then,
 any solutions of $P_{\kappa,j}(T)=0$ are contained in
 $k(\etabar(i))[\![z_d,w_d]\!]+k(\etabar(i))[\![z]\!](\!(w_d)\!)$.
\end{itemize}
We remark that the completion of
$\pi_{i,\etabar(i)}^{\ast}(\nbige,\nabla)$ at such $Q$
is mildly ramified,
which can be shown by the same argument as the proof of
Lemma \ref{lem;07.4.14.20}.

\vspace{.1in}

Due to a theorem of Sabbah in \cite{sabbah4},
we can take a regular birational map
$F:\Xbar\lrarr X^{(i_1)}$ as follows:
\begin{itemize}
\item
$F$ is the blow up
along the ideal supported at the cross points
of the divisor $\pi_{i_1}^{-1}(P)$.
\item
Any cross points of the divisor $G^{-1}(P)$
are not turning points for
$(\nbigebar,\nablabar):=G^{\ast}(\nbige,\nabla)$,
where $G:=\pi_{i_1}\circ F$.
\end{itemize}
Let $Q$ be a point of the smooth part of
$G^{-1}(P)\subset \Xbar$
which is a turning point for $(\nbigebar,\nablabar)$.
We remark that $F(Q)\in X^{(i_1)}$
is contained in some exceptional divisor
with respect to $\pi^{(i_1)}$.
We take a subring $R_0\subset\cnum$
finitely generated over $R^{(i_1)}$,
on which $Q$ is defined.
We may also have a neighbourhood $U_Q$
with an \'{e}tale map $(u,v):U_Q\lrarr A^2$
around $Q$
such that $v^{-1}(0)=G^{-1}(P)\cap U_Q$.
By considering the completion at $Q$,
we obtain the free $R_0[\![u]\!](\!(v)\!)$-module
$\nbigehat_{R_0}$
with a meromorphic connection $\nablahat_{R_0}$.

\begin{lem}
\label{lem;07.4.15.3}
$(\nbigehat_{R_0},\nablahat_{R_0})$ is mildly ramified.
\end{lem}
\pf
Let $\etabar_0$ be a geometric point of $\Spec R_0$
over some $\etabar(i_1)\in S(R^{(i_1)})$. 
We have only to show that 
$(\nbigehat_{\etabar_0},\nablahat_{\etabar_0})$
is mildly ramified.
Assume $F(Q)$ is a cross point of
the divisor $\pi_{i_1}^{-1}(P)$.
Then, $F_{\etabar_0}(Q_{\etabar_0})$ is not contained
in any $Z^{(i_1)}_{\kappa,j}$,
and hence the ramification around $Q_{\etabar_0}$ may occur
only at $G_{\etabar_0}^{-1}(P_{\etabar_0})$.
In the case where $F(Q)$ is contained in the smooth part of
$\pi_{i_1}^{-1}(P)$,
the claim follows from our choice of $i_1$.
\hfill\qed

\vspace{.1in}
Then, we can control the irregular values
for $(\nbigebar,\nablabar)$.
\begin{lem}
\label{lem;07.4.28.10}
\mbox{{}}
Let $S$ denote the set of the irregular values of
$(\nbigebar,\nablabar_v)\otimes\cnum(\!(u)\!)(\!(v)\!)$.
\begin{itemize}
\item
$S$ is contained in $\cnum[\![u]\!](\!(v_d)\!)/\cnum[\![u,v_d]\!]$
for some appropriate $d$.
\item
 For any curve $\varphi:C\lrarr \Xtilde$
 such that $\varphi(C)\cap D_Q=\{Q\}$,
 where $D_Q$ denotes the exceptional divisor containing $Q$,
 the irregular values of 
 $\varphi^{\ast}(\nbigebar,\nablabar)$
 are given by the negative parts of
 $\varphi^{\ast}\gminia$
 $(\gminia\in S)$.
\end{itemize}
\end{lem}
\pf
It follows from Proposition \ref{prop;07.4.14.15},
Proposition \ref{prop;07.6.19.1}
and Lemma \ref{lem;07.4.15.3}.
\hfill\qed

\vspace{.1in}

Now, we use the classical topology.
Let $\nbigu$ be a neighbourhood of $Q$ in $\overline{X}$.
We will shrink $\nbigu$ 
without mention in the following argument,
if it is necessary.
Let $\varphi:\nbigutilde\lrarr \nbigu$ be the ramified covering
given by $(u,v_d)\longmapsto(u,v_d^d)$ for some appropriate $d$.
We put $\nbigg:=\seisuu/d\seisuu$ which naturally acts $\nbigutilde$.
We put $\nbigd_d:=\{v_d=0\}$.
Let $M(\nbigutilde)$ (resp. $H(\nbigutilde)$) denote the space of
meromorphic (resp. holomorphic) functions
whose poles are contained in $\nbigd_d$.
For each $\gminia\in M(\nbigutilde)/H(\nbigutilde)$,
we use the same notation to denote 
the natural lift to $M(\nbigutilde)_{<0}$.
Because of Lemma \ref{lem;07.4.28.10},
there exists the finite subset $S\subset M(\nbigu_d)/H(\nbigu_d)$
which gives the irregular values of
$(\nbigebar,\nablabar_v)\otimes\cnum(\!(u)\!)(\!(v)\!)$.
(Meromorphic property of the irregular values
 is shown in Theorem 2.3.1 of \cite{sabbah4}, for example.)
Let $S_1$ denote the set of pairs
$(\gminia,\gminib)\in S^2$ such that $\gminia\neq \gminib$.

We put $\varphi(\gminia):=\prod_{\sigma\in \nbigg}
 \sigma^{\ast}\gminia$ for any $\gminia\in S$
which give the meromorphic functions $\varphi(\gminia)$ on $\nbigu$.
For any $(\gminia,\gminib)\in S_1$,
we have the meromorphic functions
$\varphi(\gminia-\gminib)$ on $\nbigu$,
similarly.
The union of the zero and the pole of $\varphi(\gminia)$
is denoted by $|\varphi(\gminia)|$.
We use the notation 
$|\varphi(\gminia-\gminib)|$
in a similar meaning.

We can take the resolution $\kappa:\nbigu_1\lrarr \nbigu$
such that the following holds:
\begin{itemize}
\item
 $\kappa^{-1}\bigl(|\varphi(\gminia)|\cup D_Q\bigr)$
 and $\kappa^{-1}\bigl(|\varphi(\gminia-\gminib)|\cup D_Q\bigr)$
 are normal crossing
 for any $\gminia\in S$ and $(\gminia,\gminib)\in S_1$.
 Here $D_Q$ denotes the component of
 $G^{-1}(P)$ such that $Q\in D_Q$.
\item
 The zero and the pole of $\kappa^{-1}(\varphi(\gminia))$
 have no intersections for any $\gminia\in S$.
 The zero and the pole of
     $\kappa^{-1}\bigl(\varphi(\gminia-\gminib)\bigr)$
 have no intersections 
 for any $(\gminia,\gminib)\in S_1$.
\item
 For any $(\gminia,\gminib),(\gminia',\gminib')\in S_1$,
 the ideals generated by
 $\kappa^{-1}(\varphi(\gminia-\gminib))$
 and $\kappa^{-1}(\varphi(\gminia'-\gminib'))$
 are principal.
\end{itemize}

Applying Sabbah's theorem,
we can take $\nu:\nbigu''\lrarr \nbigu'$ such that
any cross points of the divisor
$(\kappa\circ\nu)^{-1}(Q)$ are not turning.
We put $\kappatilde:=\kappa\circ\nu$,
for which the above three conditions are satisfied.
For any point $Q'$ of the smooth part of $\kappatilde^{-1}(Q)$,
the irregular values of
$\kappatilde^{-1}(\nbigebar,\nablabar)$ around $Q'$
are given by the negative parts of
$\kappatilde^{-1}(\gminia)$
due to Lemma \ref{lem;07.4.28.10}.
By using Proposition \ref{prop;07.9.20.1},
we can conclude that 
$Q'$ is not a turning point.
Therefore, we have no turning points in
$\kappatilde^{-1}(Q)$.
Applying the procedure to 
any turning points for $(\nbigebar,\nablabar)$
contained in $G^{-1}(U)$,
we can resolve them.
Thus the proof of Theorem \ref{thm;07.4.15.10}
is finished.
\hfill\qed

\noindent
{Address\\
Department of Mathematics,
Kyoto University,
Kyoto 606-8502, Japan,\\
takuro@math.kyoto-u.ac.jp
}


\begin{thebibliography}{99}


\bibitem{andre}
Y. Andr\'{e},
 {\em Structure des connexions m\'eromorphes formelles
 de plusieurs variables et semi-continuit\'e de l'irr\'egularit\'e},
 AG/0701895

\bibitem{artin}
E. Artin,
{\em Algebraic Numbers and Algebraic Functions},
 Gordon and Breach, Science Publishers,
 1967

\bibitem{bingener}
J. Bingener,
{\em \"Uber Formale Komplexe R\"aume},
Manuscripta Math. {\bf 24}, (1978),
253--293.

\bibitem{d}
P. Deligne, {\it
\'{E}quations diff\'{e}rentielles \`{a} points 
 singuliers r\'{e}guliers},
Lectures Notes  in Mathematics, {\bf 163}, Springer-Verlag,
Berlin-Heidelberg-New York, 1970.

\bibitem{grothendieck-murre}
A. Grothendieck and J. Murre,
{\em The tame fundamental group
 of a formal neighbourhood of a divisor
 with normal crossings on a scheme},
LNM {\bf 208}, Springer-Verlag,
Berlin-Heidelberg-New York,
1971


\bibitem{har1}
R. Hartshorne,
{\em Algebraic geometry},
Springer-Verlag, New York-Heidelberg, 1977

\bibitem{hukuhara}
M. Hukuhara,
{\em Sur les points singuliers des \'equations diff\'erentielles
 lin\'eaires, II},
 Jour. Fac. Sci. Hokkaido Univ.,
 {\bf 5}, (1937), 123--166.


\bibitem{k5}
M. Kashiwara,
{\em Semisimple holonomic $D$-modules},
in {\em Topological Field Theory, Primitive Forms
and Related Topics}, Progress in Math, vol {\bf 160},
Birkh\"{a}user, %\"
(1998), 267--271.


\bibitem{katz1}
N. Katz,
{\em Nilpotent connections and the monodromy theorem;
 applications of a result of Turrittin},
 Publ. Math. I.H.E.S. {\bf  39}, (1970), 175--232.

\bibitem{katz2}
N. Katz,
{\em Algebraic solutions of differential equations
 ($p$-curvature and the Hodge filtration)},
 Invent. math. {\bf 18} (1972), 1--118.

\bibitem{katz3}
N. Katz,
{\em A conjecture in the arithmetic theory
 of differential equations},
 Bulletin de S.M.F. {\bf 110}, (1982),
203--239

\bibitem{laszlo-pauly}
Y. Laszlo and C. Pauly,
{\em On the Hitchin morphism in positive characteristic},
Intenat. Math. Res. Notices,
(2001), 129--143,
math.AG/0005044.

\bibitem{levelt}
A. Levelt, {\em Jordan decomosition
for a class of singular differential operators},
Ark. Math. {\bf 13}, (1975), 1--27

\bibitem{majima}
H. Majima,
{\em Asymptotic analysis for integrable connections
with irregular singular points},
Lecture Notes in Mathematics, {\bf 1075},
Springer-Verlag, Berlin, 1984

\bibitem{malgrange}
B. Malgrange,
{\em Connexions m\'eromorphies $2$, Le r\'eseau canonique},
 Invent. Math. {\bf 124}, (1996) 367--387.

\bibitem{milne}
J. Milne,
{\em Etale cohomology},
Princeton University Press,
1980

\bibitem{mochi6}
T. Mochizuki,
{\em Wild harmonic bundles
 and wild pure twistor $D$-modules},
in preparation.

\bibitem{mumford}
D. Mumford,
{\em Lectures on curves on an algebraic surface.}
Annals of Mathematics Studies, {\bf 59} 
Princeton University Press, Princeton, N.J. 
1966

\bibitem{ogus-vologodsky}
A. Ogus and V. Vologodsky,
{\em Nonabelian Hodge Theory in Characteristic $p$},
AG/0507476

\bibitem{robba-christol}
P. Robba and G. Christol,
{\em \'{E}quations Diff\'erentielles $p$-adiques},
Hermann, Paris, 1994.

\bibitem{sabbah4}
C. Sabbah,
{\em \'{E}quations diff\'erentielles
\`a points singuliers irr\'eguliers
et ph\'enom\`ene de Stokes
en dimension $2$},
Ast\'{e}risque, {\bf 263},
Soci\'{e}t\'{e} Math\'{e}matique
de France, Paris, 2000.

\bibitem{tsuchimoto}
Y. Tsuchimoto,
{\em Endomorphisms of Weyl algebra
and $p$-curvatures},
Osaka J. Math.,
{\bf 42},
(2005), 435--452.

\bibitem{turrittin}
H. L. Turrittin,
{\em Convergent solutions of ordinary
 homogeneous differential equations
 in the neighborhood of a singular point},
 Acta Math., {\bf 93}, (1955), 27--66.

\bibitem{wasow}
W. Wasow, 
{\em Asymptotic expansions for ordinary equations},
Reprint of 1976 edition.
Dover Publications, Inc.,
New York, 1987

\bibitem{raynaud}
{\em Rev\^{e}tements \'etales et Groupe Fondamental,
S\'eminaire de G\'eom\'etrie Alg\'ebrique du Bois Marie
 {\rm 1960--1961} SGA $1$},
Dirig\'e par Alexander Grothendieck,
Augment\'e de deux expos\'es de M. Raynaud,
Lecture Notes in Mathematics {\bf 224},
Springer-Verlag, Berlin-Heidelberg-New York,
(1971)


\end{thebibliography}
\end{document}